\documentclass[12pt]{amsart}
\usepackage{amssymb,amsxtra}
\theoremstyle{plain}

\newcommand{\cleqn}{\setcounter{equation}{0}}
\newcommand{\clth}{\setcounter{theorem}{0}}
\newcommand {\sectionnew}[1]{\section{#1}\cleqn\clth}

\newtheorem{theorem}{Theorem}[section]
\newtheorem{procl}[theorem]{}

\theoremstyle{definition}
\newtheorem{subsec}[theorem]{}

\def\L{{\mathcal L}}
\def\N{{\mathcal N}}
\def\O{{\mathcal O}}

\def\NN{{\mathbb N}}
\def\ZZ{{\mathbb Z}}

\def\ihat{\,\widehat{i}\,}
\def\jhat{\,\widehat{j}\,}
\def\lhat{\,\widehat{l}\,}
\def\qhat{\widehat{q}\,}
\def\shat{\widehat{s}\,}

\def\bfr{\mathbf{r}}

\def\itil{\tilde{i}}
\def\jtil{\tilde{j}}

\def\ltil{\tilde{l}}  

\def\xtil{\widetilde{x}}
\def\ytil{\widetilde{y}}

\def\Itil{\widetilde{I}}
\def\Itl{\tilde{I}}
\def\Jtil{\widetilde{J}}
\def\Jtl{\tilde{J}}
\def\Mtil{\widetilde{M}}
\def\Mtl{\tilde{M}}
\def\Ntil{\widetilde{N}}
\def\Ntl{\tilde{N}}

\def\Stil{\widetilde{S}}
\def\Stl{\tilde{S}}
\def\Ttil{\widetilde{T}}
\def\Ttl{\tilde{T}}
\def\Util{\widetilde{U}}

\def\lambdatil{\widetilde{\lambda}}
\def\mutil{\widetilde{\mu}}

\def\abar{\overline{a}}
\def\bbar{\overline{b}}

\def\Oq{{\mathcal O}_q}
\def\OqMn{\Oq(M_n(k))}

\def\pmqdot{\pm q^{\bullet}}
\def\kx{k^\times}
\def\id{\operatorname{id}}
\def\OMn{\O(M_n(k))}
\def\Uqsln{U_q({\mathfrak{sl}}_n(k))}
\def\OqSLn{\Oq(SL_n(k))}

\def\mdd{|}
\def\sqp{{\sqcup}}
\def\cupt{{\cup}}
\def\capt{{\cap}}
\def\setmin{{\setminus}}
\def\dlt{{\Delta}}

\def\lsind#1#2{\{{<}#1{\Vert}#2\}}
\def\grind#1#2{\{{>}#1{\Vert}#2\}}
\def\lsqind#1#2{\{{\le}#1{\Vert}#2\}}
\def\grqind#1#2{\{{\ge}#1{\Vert}#2\}}
\def\nat{^{\natural}}

\begin{document}


\title[Commutation relations for arbitrary quantum minors]
{Commutation relations for arbitrary quantum minors}

\author[K. R. Goodearl]{K. R. Goodearl}

\address{Department of Mathematics, University of California,
Santa Barbara, CA 93106, USA}

\email{goodearl@math.ucsb.edu}


\thanks{This research was partially supported by a grant from the
National Science Foundation.}

\begin{abstract} 
Complete sets of commutation relations for arbitrary pairs of
quantum minors are computed, with explicit coefficients in closed form. 
\end{abstract}

\maketitle


\section*{Introduction}

The title of this paper begins with what may seem a misnomer -- the term
{\it commutation relation\/}, in current usage, does not refer to a 
commutativity condition, $xy=yx$, but has evolved to encompass various
``skew commutativity'' conditions that have proved to be useful
replacements for commutativity. Older types of commutation relations
include conditions of the form $xy-yx=z$, used in defining Weyl algebras
and enveloping algebras. In quantized versions of classical algebras,
relations such as $xy=qyx$ (known as {\it $q$-commutation\/}) appear, along
with mixtures of both types. Thus, it has become common to refer to
any equation of the form $xy= \lambda yx +z$, where $\lambda$ is a nonzero
scalar, as a {\it commutation relation for $x$ and $y$\/}. One important
use of such a relation, especially in enveloping algebras, is that if the
algebra supports a filtration such that $\deg(z)< \deg(x)+\deg(y)$, then
the images of $x$ and $y$ in the associated graded algebra, call them
$\xtil$ and $\ytil$, commute up to a scalar: $\xtil\ytil=
\lambda\ytil\xtil$. Similarly, the cosets of $x$ and $y$ modulo the ideal
generated by $z$ commute up to $\lambda$. Such coset relations are key
ingredients in the work of Soibelman \cite{Soi}, Hodges-Levasseur
\cite{HoLethr, HoLeSLn}, Joseph \cite{Jos}, and others on quantized
coordinate rings.

In many quantized algebras, the available commutation relations are
homogeneous and quadratic, of the form $xy= \lambda yx+ \sum_i \mu_i x_iy_i$
(where $\lambda$ and the $\mu_i$ are nonzero scalars). Relations of
this type are particularly important in establishing a (noncommutative)
standard basis of monomials in generators that include the elements $x$,
$y$, $x_i$, $y_i$. Namely, if the generators are ordered in such a way that
each $x_i\le y_i$ but $x>y$, then the given relation allows one to
rewrite monomials involving $xy$ as linear combinations of monomials
closer to standard form. For example, noncommutative standard bases
have been constructed by Lakshmibai and Reshetikhin \cite{LaRe1,
LaRe2} for quantized coordinate rings of flag varieties and Schubert
schemes, by the author and Lenagan \cite{GLduke} for quantum matrix
algebras, and by Lenagan and Rigal \cite{LeRi} for quantum Grassmannians
and quantum determinantal rings.

In order to work effectively with quantized coordinate rings of
matrices, Grassmannians, special or general linear groups, and related
algebras, one needs explicit commutation relations for quantum minors and
related elements. Such relations have often been derived for special cases
as needed, either by induction on the size of the minors, using quantum
Laplace relations, as in Parshall-Wang \cite{PaWa} and Taft-Towber
\cite{TaTo}, or by applying the quasitriangular structure of
$\Uqsln$ (that is, its universal R-matrix) to coordinate functions in
$\OqSLn$, as in the work of Lakshmibai-Reshetikhin \cite{LaRe1, LaRe2},
Soibelman \cite{Soi}, and Hodges-Levasseur
\cite{HoLethr, HoLeSLn}. Along the former line, the most complete
results to date were obtained by Fioresi \cite{Fi, Fitwo}, who
developed an algorithm which yields a commutation relation for any pair
of quantum minors. This algorithm is an iterative procedure, in which
certain products of quantum minors may appear multiple times; explicit
coefficients are produced, but are not expressed as closed formulas.
Via the quasitriangular approach, general commutation relations for pairs of
coordinate functions in quantized coordinate rings $\Oq(G)$, where $G$ is a
semisimple Lie group, have been derived in special cases (e.g., see
\cite{LaRe1, LaRe2, Soi, HoLethr, HoLeSLn}), not all with explicit
coefficients. (Quantum minors in $\OqSLn$ are special coordinate functions.)
Perhaps the largest group of explicit commutation relations obtained in
this way appeared in Hodges-Levasseur-Toro \cite{HLT} (cf.~also
\cite{BrGo}). However, to make these fully explicit, canonical elements for
the Rosso-Tanisaki Killing form on $\Uqsln$ had to be computed.

Here we introduce a new method -- new only in the sense that it has
apparently not been used for this purpose before -- with which we derive
complete commutation relations for arbitrary pairs of quantum minors,
with explicit coefficients in closed form. Our method is dual to the
quasitriangular approach, as it relies on the coquasitriangular (or
braided) bialgebra structure on the quantized coordinate ring of
$n\times n$ matrices. Representation-theoretically, the two approaches
are based on equivalent information, in that a quasitriangular
(respectively, coquasitriangular) structure on a bialgebra encodes
braiding isomorphisms $V\otimes W \xrightarrow{\,\cong\,} W\otimes V$
for  finite dimensional modules (respectively, comodules) $V$ and $W$.
To record such isomorphisms, one typically requires formulas for
matrix entries. However, in the case of a coquasitriangular bialgebra
$A$, the above isomorphism information is stored more compactly, in a
bilinear form
$\bfr$ on $A$ -- the braiding isomorphism for left $A$-comodules $V$
and $W$ is then given by the formula
\[
v\otimes w \longmapsto \sum_{(v),(w)} \bfr(v_0,w_0) w_1\otimes v_1,
\]
where we have used the Sweedler notation $v\mapsto \sum_{(v)} v_0\otimes
v_1$ for the comodule structure map $V\rightarrow A\otimes V$, and similarly
for
$W$. The resulting commutation relations are equations with values of
$\bfr$ as coefficients, namely
\begin{equation} \label{eq0.1}
\sum_{(a),(b)} \bfr(a_1,b_1)a_2b_2 = \sum_{(a),(b)}
\bfr(a_2,b_2)b_1a_1 
\end{equation}
for $a,b\in A$, where now the Sweedler notation is used for the
comultiplication map $A\rightarrow A\otimes A$.

When $A$ is the bialgebra $\OqMn$ and $a= [I\mdd J]$ and $b= [M\mdd N]$
are quantum minors (see below for notation), equation \eqref{eq0.1} becomes
\begin{equation} \label{eq0.2}
\sum_{\substack{|S|=|I|\\ |T|=|M|}} \bfr\bigl( [I\mdd S], [M\mdd T]
\bigr) [S\mdd J][T\mdd N]=  \sum_{\substack{|S|=|J|\\ |T|=|N|}}
\bfr\bigl(  [S\mdd J], [T\mdd N] \bigr) [M\mdd T][I\mdd S]. 
\end{equation}
Observe that $[I\mdd J][M\mdd N]$ occurs on the left hand
side of \eqref{eq0.2} when $S=I$ and $T=M$, while $[M\mdd N][I\mdd J]$
occurs on the right when $S=J$ and $T=N$. As we shall see, the coefficients
for these terms, namely $\bfr\bigl( [I\mdd I], [M\mdd M] \bigr)$ and
$\bfr\bigl( [J\mdd J], [N\mdd N] \bigr)$, are nonzero (in fact, they
are powers of $q$). Thus, to obtain explicit commutation relations for
$[I\mdd J]$ and
$[M\mdd N]$, we only need to compute the values $\bfr\bigl( [I\mdd S],
[M\mdd T]
\bigr)$ and $\bfr\bigl(  [S\mdd J], [T\mdd N] \bigr)$. This is
precisely what we do in the paper -- see especially Theorems \ref{thm4.6}
and \ref{thm5.2}. Additional relations follow from these by various
symmetries, or by investing quantum Laplace relations. (Quantum Pl\"ucker
relations in quantum Grassmannians can also be used for this purpose.) See
Theorems \ref{thm5.6}, \ref{thm6.3} and Corollaries \ref{cor5.4},
\ref{cor5.7}, \ref{cor6.4}.
\medskip

Our notation and conventions are collected in Section \ref{sec1}. In
particular, the relations we use for $\OqMn$ are displayed in \eqref{eq1.6},
so that the reader may compare with other papers in which $q$ is replaced
by $q^{-1}$ or $q^2$. Our computations of the values of the form $\bfr$
on pairs of quantum minors occupy Sections \ref{sec2} and \ref{sec4}; the
intermediate Section \ref{sec3} provides a first set of commutation
relations to illustrate our methods. The general commutation relations are
derived in Sections \ref{sec5} and \ref{sec6}, and we conclude by using
these relations, in Section \ref{sec7}, to evaluate the standard Poisson
bracket on pairs of classical minors.


\sectionnew{Notation and conventions} \label{sec1}

Fix a positive integer $n$, a base field $k$, and a nonzero scalar
$q\in \kx$. We work within the standard single-parameter quantized
coordinate ring of $n\times n$ matrices over $k$, which we denote
$\OqMn$, as defined in \S\ref{sub1.2} below. We use the abbreviation
\begin{equation} \label{eq1.1}
\qhat= q-q^{-1}, 
\end{equation}
since this scalar appears in numerous formulas.

\begin{subsec} \label{sub1.1} \textbf{R-matrix.} 
The standard R-matrix of type
$A_{n-1}$ can be presented in the form
\begin{equation} \label{eq1.2}
R= q\sum_{i=1}^n e_{ii}\otimes e_{ii}+ \sum^n_{\substack{i,j=1\\ i\ne j}}
e_{ii}\otimes e_{jj} +\qhat \sum^n_{\substack{i,j=1\\ i> j}} e_{ij}\otimes
e_{ji} 
\end{equation}
\cite[Equation (1.5), p.~200]{RTF}. We view $R$ as a linear
automorphism of $k^n\otimes k^n$, which acts on the standard basis
vectors $x_i\otimes x_j$ according to the following formula, using the
conventions of \cite{KlSc}:
\begin{equation} \label{eq1.3}
R(x_l\otimes x_m)= \sum_{i,j=1}^n R^{ij}_{lm} x_i\otimes x_j. 
\end{equation} 
The entries of the $n^2\times n^2$ matrix $R^{ij}_{lm}$ are as
follows (cf.~\cite[Equation (9.13), p.~309]{KlSc}):
\begin{equation} \label{eq1.4}
\begin{aligned} R_{ii}^{ii} &= q &&\qquad (\text{all\ } i)
&\qquad\qquad\qquad R_{ij}^{ij} &= 1 &&\qquad (i\ne j) \\
R^{ij}_{ji} &= \qhat &&\qquad (i>j) &\qquad\qquad\qquad R^{ij}_{lm} &= 0
&&\qquad \text{(otherwise)}. 
\end{aligned}  
\end{equation}
\end{subsec}

\begin{subsec} \label{sub1.2} \textbf{Generators, relations, and grading.} 
The algebra $A=
\OqMn$ is obtained from \eqref{eq1.4} by the
Faddeev-Reshetikhin-Takhtadzhyan construction, namely as the
$k$-al\-ge\-bra
$A(R)$ presented by generators $X_{ij}$ (for $i,j=1,\dots,n$) and
relations
\begin{equation} \label{eq1.5}
\sum_{s,t=1}^n R^{ij}_{st} X_{sl}X_{tm}= \sum_{s,t=1}^n X_{jt}X_{is}
R^{st}_{lm}  
\end{equation}
for all $i,j,l,m= 1,\dots,n$. (See \cite[Definition 1, p.~197]{RTF}
and \cite[\S9.1.1]{KlSc}. We have written
$X_{ij}$ for the generators labelled $t_{ij}$ in \cite{RTF} and
$u^i_j$ in \cite{KlSc}.) As is well known, the relations \eqref{eq1.5} are
equivalent to
\begin{equation} \label{eq1.6}
\begin{aligned} X_{ij}X_{lj} &= qX_{lj}X_{ij} &&\qquad (i<l) \\
X_{ij}X_{im} &= qX_{im}X_{ij} &&\qquad (j<m) \\
X_{ij}X_{lm} &= X_{lm}X_{ij} &&\qquad (i<l,\ j>m) \\
X_{ij}X_{lm} - X_{lm}X_{ij} &= \qhat X_{im} X_{lj} &&\qquad (i<l,\ j<m)
\end{aligned}  
\end{equation} 
(cf.~\cite[Equations (9.17), p.~310]{KlSc}). Some authors define
quantum matrices using relations as in \eqref{eq1.6} but with $q$ replaced
by $q^{-1}$; thus, the algebras they define match what we would label
$\O_{q^{-1}}(M_n(k))$. See \cite[p.~3317]{LaTo} or \cite[Equations (3.5a),
p.~37]{PaWa}, for example. In comparing our work with those papers, we must
be careful to interchange $q$ and $q^{-1}$. However,
$\qhat$ is defined to be $q^{-1}-q$ in \cite[p.~38]{PaWa}, and so we
do not change
$\qhat$ when carrying over results from that paper.

Because of the homogeneity of the relations \eqref{eq1.6}, $A$ carries a
natural $(\ZZ^n\times\ZZ^n)$-grading, such that each $X_{ij}$ is
homogeneous of degree $(\epsilon_i,\epsilon_j)$, where
$\epsilon_1,\dots,\varepsilon_n$ are the standard basis elements for
$\ZZ^n$.
\end{subsec}

\begin{subsec} \label{sub1.3} \textbf{Coquasitriangular structure.} 
We follow \cite[Section 1]{Hay} in defining a {\it
coquasitriangular bialgebra\/} (also called a {\it bialgebra with
braiding structure\/}
\cite[Theorem 2.7]{LaTo} or a {\it cobraided bialgebra\/}
\cite[Definition VIII.5.1]{Kas}) to be a bialgebra $B$ equipped with a
convolution-invertible bilinear form $\bfr: B\otimes B\rightarrow k$
such that
\begin{align}
\sum_{(a),(b)} \bfr(a_1,b_1)a_2b_2 &= \sum_{(a),(b)}
\bfr(a_2,b_2)b_1a_1   \tag{1.7)(\text{i}} \label{eq1.7i}\\ 
\bfr(ab,c) &= \sum_{(c)} \bfr(a,c_1)\bfr(b,c_2) \tag{1.7)(\text{ii}}
\label{eq1.7ii}\\
\bfr(a,bc) &= \sum_{(a)} \bfr(a_1,c)\bfr(a_2,b) \tag{1.7)(\text{iii}}
\label{eq1.7iii}\\
\bfr(a,1) &= \bfr(1,a)= \varepsilon(a) \tag{1.7)(\text{iv}} \label{eq1.7iv}
\end{align}
for all $a,b,c\in B$, where we have written $\bfr(x\otimes y)$ as
$\bfr(x,y)$ for convenience, and have used the Sweedler notation for
comultiplication in the form
$\Delta(x)=
\sum_{(x)} x_1\otimes x_2$. Condition \eqref{eq1.7iv} is redundant by
\cite[Proposition 10.2(ii), p.~333]{KlSc}. Thus, the above definition
agrees with \cite[Definition VIII.5.1]{Kas}, \cite[Definition
10.1, pp.~331-2]{KlSc}, and
\cite[Definition 7.3.1]{LaRa}, but not with the conditions in
\cite[Theorem 2.7]{LaTo}. However, the latter conditions match those
of \eqref{eq1.7i}--(iv) if one uses the form
$\langle-|-\rangle$ given by $\langle a|b\rangle= \bfr(b,a)$.

\setcounter{equation}{7}

By \cite[Theorem 10.7, p.~337]{KlSc}, whenever $R$ is an invertible
R-matrix satisfying the original form of the quantum Yang-Baxter
equation, the FRT-algebra $A(R)$ is coquasitriangular with respect to
the form $\bfr$ determined by
\begin{equation} \label{eq1.8}
\bfr(X_{ij}, X_{lm})= R^{il}_{jm}  
\end{equation}
for all $i,j,l,m$. (By the {\it original QYBE\/}, we mean the equation
$R_{12}R_{13}R_{23}= R_{23}R_{13}R_{12}$ \cite[Equation (0.7),
p.~195]{RTF}, as opposed to the
form exhibiting the braid relation, namely
$R_{12}R_{23}R_{12}= R_{23}R_{12}R_{23}$.) Note that, in view of
\eqref{eq1.8}, if we put $a= X_{il}$ and $b= X_{jm}$ into \eqref{eq1.7i},
we recover the relations \eqref{eq1.5}.

It is well known that the R-matrix given in \eqref{eq1.2}
satisfies the original QYBE (e.g., \cite[\S8.1.2, pp.~246-7
and Equation (8.60), p.~270]{KlSc}).  Consequently:
\end{subsec}

\begin{procl} \label{thm1.4} {\bf Theorem.} 
The algebra $A= \OqMn$ is a
coquasitriangular bialgebra with respect to the bilinear form
$\bfr:A\otimes A\rightarrow k$ determined by the following conditions:
\begin{equation} \label{eq1.9}
\begin{aligned}
 \bfr(X_{ii},X_{ii}) &= q &&\qquad (\textup{all\ }i) 
&\qquad\quad \bfr(X_{ii},X_{jj}) &= 1 &&\qquad (i\ne j) \\
\bfr(X_{ij},X_{ji}) &= \qhat &&\qquad (i>j) &\qquad\quad
\bfr(X_{ij},X_{lm}) &= 0 &&\qquad \textup{(otherwise)}. \quad\square 
\end{aligned}  
\end{equation}
\end{procl}

\begin{subsec} \label{sub1.5} \textbf{Quantum minors.} 
We write $[I\mdd J]$ for the quantum
minor in
$A$ with row index set $I$ and column index set $J$; this minor is just
the quantum determinant in the subalgebra $k\langle X_{ij} \mid i\in
I,\, j\in J\rangle$, which is naturally isomorphic to $\Oq(M_{|I|}(k))$.
Specifically, if we write the elements of $I$ and $J$ in ascending
order, say 
\[
I = \{i_1< \cdots< i_t\} \qquad\qquad\qquad\qquad J = \{j_1< \cdots<
j_t\},
\]
then
\begin{equation} \label{eq1.10}
\begin{aligned}
 {[I\mdd J]} &= \sum_{\sigma\in S_t} (-q)^{\ell(\sigma)}
X_{i_{\sigma(1)},j_1} X_{i_{\sigma(2)},j_2} \cdots
X_{i_{\sigma(t)},j_t} \\
 &= \sum_{\sigma\in S_t} (-q)^{\ell(\sigma)}
X_{i_1,j_{\sigma(1)}} X_{i_2,j_{\sigma(2)}} \cdots
X_{i_t,j_{\sigma(t)}}, 
\end{aligned}  
\end{equation}
where $\ell(\sigma)$ denotes the
length of the permutation $\sigma\in S_t$ as a product of simple
transpositions $(l,l+1)$ (cf.~\cite[equations (9.18) and
(9.20), pp.~311-312]{KlSc}, \cite[p.~43]{PaWa}). Note that
$[I\mdd J]$ is homogeneous of degree
\[
(\epsilon_{i_1}+ \cdots+ \epsilon_{i_t}, \epsilon_{j_1} +\cdots+
\epsilon_{j_t})
\]
with respect to the grading of \S\ref{sub1.2}.

Comultiplication of quantum minors is given by the rule
\begin{equation} \label{eq1.11}
\Delta \bigl( [I\mdd J] \bigr)= \sum_{\substack{K\subseteq \{1,\dots,n\}\\
|K|=|I|}} [I\mdd K] [K\mdd J]  
\end{equation}
(e.g., \cite[Proposition 9.7(ii), p.~312]{KlSc}).
\end{subsec}

\begin{subsec} \label{sub1.6} \textbf{Transpose and anti-transpose.} 
As observed in
\cite[Proposition 3.7.1(1)]{PaWa}, there is a $k$-algebra automorphism
$\tau$ on
$A$ such that
$\tau(X_{ij})= X_{ji}$ for all $i,j$. We refer to $\tau$ as the {\it
transpose automorphism\/}. There is also a $k$-algebra
anti-automorphism $\tau_2$ on $A$ sending $X_{ij} \mapsto
X_{n+1-i,n+1-j}$ for all $i,j$ \cite[Proposition 3.7.1(2)]{PaWa}. This
proposition also shows that $\tau$ is a coalgebra anti-automorphism
while $\tau_2$ is a coalgebra automorphism, that is,
\begin{alignat}{2}
 \Delta\circ\tau &= \phi\circ(\tau\otimes\tau)\circ\Delta \qquad\qquad\qquad
&\Delta\circ\tau_2 &= (\tau_2\otimes\tau_2)\circ \Delta,
\notag \end{alignat} 
where $\phi$ is the {\it flip\/} automorphism on
$A\otimes A$, sending
$a\otimes b\mapsto b\otimes a$ for all $a,b\in A$. Hence, 
\begin{alignat}{2}
\Delta\tau(a) &= \sum_{(a)} \tau(a_2)\otimes \tau(a_1) \qquad\qquad\qquad 
&\Delta\tau_2(a) &= \sum_{(a)} \tau_2(a_1)\otimes \tau_2(a_2)
\notag \end{alignat} 
for $a\in A$.
Consequently, when writing out $\Delta\tau(a)$ and $\Delta\tau_2(a)$ in
Sweedler notation we may take
\begin{equation}\label{eq1.12}
\begin{aligned}
  \tau(a)_1 &= \tau(a_2) &\qquad\qquad\qquad \tau(a)_2 &=
\tau(a_1) \\
\tau_2(a)_1 &= \tau_2(a_1) &\qquad\qquad\qquad\qquad \tau_2(a)_2 &=
\tau_2(a_2).
\end{aligned}   
\end{equation}
We recall from \cite[Lemma 4.3.1]{PaWa} that 
\begin{equation} \label{eq1.13}
 \tau\bigl( [I\mdd J] \bigr) = [J\mdd I] 
\qquad\qquad\qquad
\tau_2\bigl(
[I\mdd J] \bigr) = [\omega_0I\mdd \omega_0J]  
\end{equation}
for all quantum minors $[I\mdd J]$ in $A$, where $\omega_0$ is the
longest element of $S_n$, that is, the permutation $i\mapsto n+1-i$.

As discussed in \cite[Remark 3.7.2]{PaWa}, there is an
isomorphism (of bialgebras) $\OqMn\rightarrow \O_{q^{-1}}(M_n(k))$ that
sends
$X_{ij}
\mapsto X'_{n+1-i,n+1-j}$ for all $i,j$, where the
$X'_{\bullet,\bullet}$ are the standard generators for
$\O_{q^{-1}}(M_n(k))$. Let us call this isomorphism $\beta$, and let us
use the notation $[I\mdd J]'$ for quantum minors in
$\O_{q^{-1}}(M_n(k))$. It was shown in \cite[proof of
Corollary 5.9]{GLjalg} that
\begin{equation} \label{eq1.14}
\beta\bigl( [I\mdd J] \bigr)= [\omega_0I\mdd \omega_0J]'   
\end{equation}
for all quantum minors $[I\mdd J]$ in $A$.
\end{subsec}

\begin{procl} \label{lem1.7} {\bf Lemma.} 
The form $\bfr$ satisfies
$\bfr(a,b)=
\bfr(\tau(b),\tau(a))= \bfr(\tau_2(b),\tau_2(a))$ for all
$a,b\in A$. In particular,
\begin{equation} \label{eq1.15}
\bfr\bigl( [I\mdd J], [M\mdd N] \bigr)= \bfr\bigl( [N\mdd M], [J\mdd
I] \bigr)= \bfr\bigl( [\omega_0M\mdd 
\omega_0N], [\omega_0I\mdd \omega_0J] \bigr)
\end{equation} 
for all quantum minors $[I\mdd J]$ and $[M\mdd N]$ in $A$.
 \end{procl}

\begin{proof} Set $\bfr'(a,b)= \bfr(\tau(b),\tau(a))$ and $\bfr''(a,b)=
\bfr(\tau_2(b),\tau_2(a))$ for all
$a,b\in A$, and note from \eqref{eq1.9} that $\bfr'(X_{ij},X_{lm})=
\bfr''(X_{ij},X_{lm})= \bfr(X_{ij},X_{lm})$ for all $i,j,l,m$. To prove
that $\bfr'$ and $\bfr''$ coincide with $\bfr$, it suffices to show that
these forms agree on all monomials in the $X_{ij}$. This will
be clear by induction on the lengths of the monomials once we show that
$\bfr'$ and $\bfr''$ satisfy \eqref{eq1.7ii} and \eqref{eq1.7iii}. These
identities are routine with the aid of \eqref{eq1.12}; we give one sample:
\begin{equation}
\begin{aligned}
\bfr'(ab,c) &= \bfr(\tau(c), \tau(a)\tau(b))= \sum_{(\tau(c))} 
\bfr(\tau(c)_1, \tau(b)) \bfr(\tau(c)_2, \tau(a)) \\
 &= \sum_{(c)} \bfr(\tau(c_2), \tau(b)) \bfr(\tau(c_1), \tau(a))=
\sum_{(c)} \bfr'(b,c_2) \bfr'(a,c_1) \\
 &= \sum_{(c)} \bfr'(a,c_1) \bfr'(b,c_2)
\notag \end{aligned} 
\end{equation}
for all $a,b,c\in A$.
\end{proof}

\begin{subsec} \label{newsub1.8} \textbf{Definition of quantities
$\ell(S;T)$.}
Many formulas concerning quantum minors require powers of $q$ or $-q$ whose
exponents are quantities which might be called the number of inversions
between two sets. We follow \cite{NYM} in defining
\begin{equation} \label{eq1.16}
\ell(S;T) = \big| \{ (s,t)\in S\times T \mid s>t \} \big| 
\end{equation}
for any subsets $S,T\subseteq \{1,\dots,n\}$.
\end{subsec}

\begin{subsec} \label{sub1.8} \textbf{Quantum Laplace relations.} 
We shall need the
following {\it $q$-Laplace relations\/} from
\cite[Proposition 1.1]{NYM}, for index sets $I,J \subseteq
\{1,\dots,n\}$ of the same cardinality. If $I_1$, $I_2$ are
nonempty subsets of $I$ with $|I_1|+|I_2| = |I|$, then
\begin{equation} \label{eq1.17}
\sum_{\substack{J=J_1\sqp J_2\\
|J_l|=|I_l|}} (-q)^{\ell(J_1;J_2)} [I_1\mdd J_1][I_2\mdd J_2] =
\begin{cases}
 (-q)^{\ell(I_1;I_2)} [I\mdd J] &(I_1\capt I_2= \varnothing) \\
0  &(I_1\capt I_2\ne \varnothing),  \end{cases}
\end{equation}
while if $J_1$, $J_2$ are
nonempty subsets of $J$ with $|J_1|+|J_2| = |J|$, then
\begin{equation} \label{eq1.18}
\sum_{\substack{I=I_1\sqp I_2\\
|I_l|=|J_l|}} (-q)^{\ell(I_1;I_2)} [I_1\mdd J_1][I_2\mdd J_2] =
\begin{cases}
 (-q)^{\ell(J_1;J_2)} [I\mdd J] &(J_1\capt J_2= \varnothing) \\
0  &(J_1\capt J_2\ne \varnothing).  \end{cases}
\end{equation}
Observe that \eqref{eq1.17} holds trivially in
case $I_1$ or $I_2$ is empty, and that \eqref{eq1.18} holds trivially in
case $J_1$ or $J_2$ is empty

Reduction formulas for values of the form $\bfr$ can be obtained by
combining \eqref{eq1.17} and \eqref{eq1.18} with \eqref{eq1.7ii}(iii). For
example, if $J= J_1\sqcup J_2$, then \eqref{eq1.18} together with
\eqref{eq1.7ii} yields
\begin{equation} \label{eq1.19}
\begin{aligned}  
(-q)^{\ell(J_1;J_2)} &\bfr\bigl( [I\mdd J], [M\mdd N]
\bigr) = \\
 &\sum_{\substack{I=I_1\sqp I_2}} \sum_L (-q)^{\ell(I_1;I_2)} \bfr\bigl(
[I_1\mdd J_1], [M\mdd L] \bigr)  \bfr\bigl( [I_2\mdd J_2], [L\mdd N]
\bigr) 
\end{aligned}  
\end{equation}
for all $[M\mdd N]$.
\end{subsec} 

\begin{subsec} \label{sub1.9} \textbf{Some further notation.} 
To simplify notation for
operations on index sets, we often omit braces from singletons -- in
particular, we write
\begin{equation} \label{eq1.20} 
I\setmin i = I\setmin\{i\}  \qquad\qquad
I\sqp l = I\sqp\{l\}   \qquad\qquad
I\setmin i\sqp l = \bigl( I\setmin\{i\} \bigr) \sqp \{l\}
\end{equation} 
for $i\in I$ and $l\notin I$.
The Kronecker
delta symbol will be applied to index sets as well as to individual
indices -- thus, $\delta(I,J)=1$ when $I=J$ while
$\delta(I,J)=0$ when $I\ne J$. In the case of an index versus an index
set, the Kronecker symbol will be used to indicate membership, that is,
$\delta(i,I)=1$ means $i\in I$ while $\delta(i,I)=0$ means $i\notin I$.

Finally, we shall need the following partial order on index sets of the
same cardinality. If $I$ and $J$ are $t$-element subsets of
$\{1,\dots,n\}$, write their elements in ascending order, say
\[
I = \{i_1< i_2< \cdots< i_t\}   \qquad\qquad\qquad
J  = \{j_1< j_2< \cdots< j_t\}, 
\]
and then define
\begin{equation} \label{eq1.21}
I \le J \qquad\iff\qquad i_l\le j_l \text{\ for\ } l=1,\dots,t. 
\end{equation}
\end{subsec}


\sectionnew{Initial computations} \label{sec2}

Throughout this section, let $i$ and $j$ denote indices in
$\{1,\dots,n\}$, and let $I$, $J$, $M$, $N$ denote index sets
contained in $\{1,\dots,n\}$, with $|I|=|J|$ and $|M|=|N|$.

\begin{procl} \label{lem2.1} {\bf Lemma.} 
$\bfr\bigl( X_{ii}, [I\mdd J] \bigr) =
\bfr\bigl( [I\mdd J],X_{ii} \bigr) = q^{\delta(i,I)}
\delta(I,J)$.
\end{procl}

\begin{proof} Write $I= \{i_1<\cdots<i_t\}$ and $J= \{j_1<\cdots<j_t\}$,
and note using \eqref{eq1.10} and \eqref{eq1.7ii} that
\begin{equation} \label{eq2.1}
\begin{aligned}
\bfr&\bigl( [I\mdd J],X_{ii} \bigr) = \\ 
 &\sum_{\sigma\in S_t}
(-q)^{\ell(\sigma)} \sum_{l_1,\dots, l_{t-1}}^n
\bfr(X_{i_1j_{\sigma(1)}},X_{il_1})
\bfr(X_{i_2j_{\sigma(2)}},X_{l_1l_2}) \cdots
\bfr(X_{i_tj_{\sigma(t)}},X_{l_{t-1}i}). 
\end{aligned}
\end{equation}
In view of \eqref{eq1.9}, a nonzero term can occur in the second summation
of \eqref{eq2.1} only when
$i\le l_1\le l_2\le
\cdots\le l_{t-1}\le i$, that is, when $l_1=\cdots =l_{t-1}=i$. Hence,
\eqref{eq2.1} reduces to
\begin{equation} \label{eq2.2}
\bfr\bigl( [I\mdd J], X_{ii} \bigr) = \sum_{\sigma\in S_t}
(-q)^{\ell(\sigma)} \bfr(X_{i_1j_{\sigma(1)}},X_{ii})
\bfr(X_{i_2j_{\sigma(2)}},X_{ii}) \cdots
\bfr(X_{i_tj_{\sigma(t)}},X_{ii}). 
\end{equation}
In \eqref{eq2.2}, a nonzero term can occur in the sum only when $i_s=
j_{\sigma(s)}$ for $s=1,\dots,t$. Since the $i_s$ and $j_s$ are
arranged in ascending order, this situation only happens when $I=J$ and
$\sigma=\id$. Thus,
$\bfr\bigl( [I\mdd J], X_{ii} \bigr) =0$ when $I\ne J$, and 
\[
\bfr\bigl( [I\mdd I], X_{ii} \bigr) = \bfr( X_{i_1i_1}, X_{ii})
\bfr( X_{i_2i_2}, X_{ii}) \cdots
\bfr( X_{i_ti_t}, X_{ii})= q^{\delta(i,I)}.
\]

The formula for $\bfr\bigl( X_{ii},[I\mdd J] \bigr)$ follows via
Lemma \ref{lem1.7}.
\end{proof}

\begin{procl} \label{lem2.2} {\bf Lemma.} 
$\bfr(X_{ij},-) \equiv 0$ when $i<j$,
and
$\bfr(-,X_{ij}) \equiv 0$ when $i>j$. \end{procl}

\begin{proof} Consider any monomial $a= X_{i(1),j(1)}
X_{i(2),j(2)} \cdots X_{i(t),j(t)} \in A$. Then by \eqref{eq1.7ii},
\[
\bfr(a,X_{ij})= \sum_{l_1,\dots, l_{t-1}}^n
\bfr(X_{i(1),j(1)},X_{il_1}) \bfr(X_{i(2),j(2)},X_{l_1l_2}) \cdots
\bfr(X_{i(t),j(t)},X_{l_{t-1}j}).
\]
If some term $\bfr(X_{i(1),j(1)},X_{il_1})
\bfr(X_{i(2),j(2)},X_{l_1l_2})
\cdots \bfr(X_{i(t),j(t)},X_{l_{t-1}j})$ does not vanish, we must have
$i\le l_1\le \cdots\le l_{t-1}\le j$. This shows that $\bfr(-,X_{ij})$
can fail to vanish only when $i\le j$. The first statement of the lemma
follows via Lemma \ref{lem1.7}.
\end{proof}

\begin{procl} \label{cor2.3} {\bf Corollary.} 
$\bfr\bigl( [I\mdd J],- \bigr)
\equiv 0$ when $I
\not\geq J$, and $\bfr\bigl( -, [I\mdd J] \bigr) \equiv 0$ when $I\not\leq
J$.
\end{procl}

\begin{proof} Write $I= \{i_1<\cdots<i_t\}$ and $J= \{j_1<\cdots<j_t\}$,
and suppose that $\bfr\bigl( [I\mdd J],c \bigr) \ne 0$ for some $c\in A$.
Then by \eqref{eq1.10} and \eqref{eq1.7ii},
\[
\sum_{(c)} \bfr(X_{i_1j_{\sigma(1)}},c_1)
\bfr(X_{i_2j_{\sigma(2)}},c_2) \cdots \bfr(X_{i_tj_{\sigma(t)}},c_t) \ne
0
\]  
for some $\sigma\in S_t$. Lemma \ref{lem2.2} then implies that
$i_s\ge j_{\sigma(s)}$ for $s=1,\dots,t$.

First, $i_1\ge j_{\sigma(1)}\ge j_1$. Now let $1< s\le t$. If
$\sigma(s)\ge s$, then $i_s\ge j_{\sigma(s)}\ge j_s$. If $\sigma(s)< s$,
then $\sigma(u)\ge s$ for some $u< s$, whence
$i_s> i_u\ge j_{\sigma(u)}\ge j_s$. Thus, $i_s\ge j_s$ for all $s$,
and therefore $I\ge J$. Similarly, if $\bfr\bigl( -,[I\mdd J] \bigr)$
does not vanish, then $I\le J$. \end{proof}

\begin{procl} \label{prop2.4} {\bf Proposition.} 
If $i<j$, then
\begin{align}
\bfr\bigl( [I\mdd J], X_{ij} \bigr) &= \qhat (-q)^{|[1,i)\capt
J|- |[1,j)\capt I|} \delta(i,J) \delta(j,I) \delta(I\setmin j,J\setmin
i)  
\tag{2.3)(i}  \label{eq2.3i}  \\
 &= \qhat (-q)^{-|(i,j)\capt I\capt J|} \delta(i,J) \delta(j,I)
\delta(I\setmin j,J\setmin i).  
\tag{2.3)(ii}  \label{eq2.3ii}
\end{align}
\end{procl}

\setcounter{equation}{3}

\begin{proof} Note first that \eqref{eq2.3ii} follows from \eqref{eq2.3i}.
For if the right hand side of \eqref{eq2.3i} is nonzero, then $I= (I\capt
J)\sqp j$ and
$J= (I\capt J)\sqp i$, whence $[1,i)\capt J= [1,i)\capt I\capt J=
[1,i]\capt I\capt J$ and $[1,j)\capt I= [1,j)\capt I\capt J$.

We induct on $|I|$, the case $|I|=1$ being clear from
\eqref{eq1.9}. Now assume that $|I|>1$, and suppose that $\bfr\bigl( [I\mdd
J], X_{ij} \bigr) \ne 0$.

Choose $s\in I$, and write $I= I_1\sqp I_2$ with $I_1=\{s\}$ and $I_2=
I\setmin \{s\}$. The
$q$-Laplace relation \eqref{eq1.17} yields
\begin{equation} \label{eq2.4}
(-q)^{|[1,s)\capt I|} [I\mdd J]= \sum_{t\in J} (-q)^{|[1,t)\capt J|}
X_{st} [I\setmin s\mdd J\setmin t]. 
\end{equation}
For each $t\in J$, we have
\begin{equation} \label{eq2.5}
\bfr\bigl( X_{st}[I\setmin s\mdd J\setmin t], X_{ij} \bigr)=
\sum_{l=1}^n \bfr(X_{st}, X_{il}) \bfr\bigl( [I\setmin s\mdd J\setmin
t], X_{lj} \bigr).
\end{equation}
Since $\bfr\bigl( [I\mdd J],
X_{ij} \bigr) \ne 0$, we must have $\bfr(X_{st}, X_{il}) \bfr\bigl(
[I\setmin s\mdd J\setmin t], X_{lj} \bigr) \ne 0$ for some $l\in
\{1,\dots,n\}$ and $t\in J$.

Suppose that $i\notin J$. Then $t\ne i$, and so because
$\bfr( X_{st}, X_{il}) \ne 0$, we must have $t=s$ and $l=i$. Then
$\bfr\bigl( [I\setmin s\mdd J\setmin s], X_{ij} \bigr) \ne 0$, which
contradicts the induction hypothesis because $i\notin J\setmin s$.
Therefore $i\in J$.

Next, suppose that $j\notin I\setmin s$. If $l<j$, we would have
$\bfr\bigl( [I\setmin s\mdd J\setmin t], X_{lj} \bigr) =0$ by the
induction hypothesis. Since $\bfr(-, X_{lj})$ would vanish if $l>j$, we
must have
$l=j$. Now $\bfr(X_{st}, X_{ij})\ne 0$, and so $s=j$ and $t=i$. Thus,
either $j\in I\setmin s$ or $j=s$, so in any case we conclude that
$j\in I$.

We may now assume that $s=j$. Since $j\notin I\setmin j$, we have
$\bfr\bigl( [I\setmin j\mdd J\setmin t], X_{ij} \bigr) =0$ for all $t\in
J$ by the induction hypothesis. On the other hand, $\bfr(X_{jt}, X_{il})
=0$ for $l\ne i,j$, and $\bfr(X_{jt}, X_{ij}) =0$ for $t\ne i$. Hence,
the right hand side of \eqref{eq2.5} vanishes when $t\ne i$, and it equals
$\qhat \bfr\bigl( [I\setmin j\mdd J\setmin i], X_{jj} \bigr)$ when $t=i$.
Combining \eqref{eq2.4} and \eqref{eq2.5} thus yields
\begin{equation} \label{eq2.6}
(-q)^{|[1,j)\capt I|} \bfr\bigl( [I\mdd J], X_{ij} \bigr)=
(-q)^{|[1,i)\capt J|}
\qhat \bfr\bigl( [I\setmin j\mdd J\setmin i], X_{jj} \bigr). 
\end{equation}
Since the left hand side of \eqref{eq2.6} is nonzero by assumption, Lemma
\ref{lem2.1} implies that $I\setmin j= J\setmin i$ and $\bfr\bigl(
[I\setmin j\mdd J\setmin i], X_{jj} \bigr) =1$. The formula \eqref{eq2.3i}
follows, and the induction step is established. \end{proof}

\begin{procl} \label{cor2.5} {\bf Corollary.} 
If $i>j$, then
\begin{align} 
\bfr\bigl(X_{ij}, [I\mdd J] \bigr) &= \qhat (-q)^{|[1,j)\capt I|-
|[1,i)\capt J|} \delta(i,J) \delta(j,I) \delta(I\setmin j,J\setmin i) 
\tag{2.7)(i} \label{eq2.7i}  \\
 &= \qhat (-q)^{-|(j,i)\capt I\capt J|} \delta(i,J) \delta(j,I)
\delta(I\setmin j,J\setmin i).  
\tag{2.7)(ii} \label{eq2.7ii}
\end{align}
\end{procl}

\setcounter{equation}{7}

\begin{proof} Apply Lemma \ref{lem1.7} to Proposition \ref{prop2.4}.
\end{proof}

\begin{procl} \label{prop2.6} {\bf Proposition.} 
$\bfr\bigl( [I\mdd I], [M\mdd N]
\bigr)=
\bfr\bigl( [M\mdd N], [I\mdd I] \bigr)= q^{|I\capt M|} \delta(M,N)$.
\end{procl}

\begin{proof} This is parallel to the proof of Lemma \ref{lem2.1}. 
Write $M= \{m_1<\cdots<m_t\}$ and $N= \{n_1<\cdots<n_t\}$,
and note that
\begin{equation} \label{eq2.8}
\bfr\bigl( [M\mdd N], [I\mdd I] \bigr) = \sum_{\sigma\in S_t}
(-q)^{\ell(\sigma)} \bfr\bigl( X_{m_1n_{\sigma(1)}}
X_{m_2n_{\sigma(2)}} \cdots X_{m_tn_{\sigma(t)}}, [I\mdd I] \bigr),
\end{equation} 
while for each $\sigma\in S_t$ we have
\begin{equation} \label{eq2.9}
\begin{aligned} 
\bfr\bigl( &X_{m_1n_{\sigma(1)}}
X_{m_2n_{\sigma(2)}} \cdots X_{m_tn_{\sigma(t)}}, [I\mdd I] \bigr) = \\
 &\sum_{L_1,\dots, L_{t-1}}
\bfr\bigl(X_{m_1n_{\sigma(1)}}, [I\mdd L_1] \bigr)
\bfr\bigl(X_{m_2n_{\sigma(2)}}, [L_1\mdd L_2] \bigr) \cdots
\bfr\bigl(X_{m_tn_{\sigma(t)}}, [L_{t-1}\mdd I] \bigr).
\end{aligned}  
\end{equation} 
Consider the right hand side of \eqref{eq2.9}. By Corollary
\ref{cor2.3}, a nonzero term can occur in that sum only when $I\le L_1\le
\cdots\le L_{t-1}\le I$, and so only when all the $L_s=I$. Thus,
\begin{equation} \label{eq2.10}
\begin{aligned}
\bfr\bigl( [M&\mdd N], [I\mdd I] \bigr) = \\
 &\sum_{\sigma\in S_t} (-q)^{\ell(\sigma)} \bfr\bigl(
X_{m_1n_{\sigma(1)}}, [I\mdd I] \bigr) \bfr\bigl( X_{m_2n_{\sigma(2)}},
[I\mdd I] \bigr) \cdots \bfr\bigl( X_{m_tn_{\sigma(t)}}, [I\mdd I]
\bigr). 
\end{aligned} 
\end{equation} 
In view of Lemma \ref{lem2.2} and Corollary \ref{cor2.5}, $\bfr\bigl(
X_{ij}, [I\mdd I] \bigr)=0$ for all $i\ne j$. Hence, a nonzero term
can occur in the right hand side of \eqref{eq2.10} only when $m_s=
n_{\sigma(s)}$ for all $s$, that is, only when $M=N$ and $\sigma=\id$.
Therefore
$\bfr\bigl( [M\mdd N], [I\mdd I] \bigr)= 0$ when $M\ne N$, while 
\[
\bfr\bigl( [M\mdd M], [I\mdd I] \bigr)= \bfr\bigl( X_{m_1m_1}, [I\mdd I]
\bigr) \bfr\bigl( X_{m_2m_2}, [I\mdd I] \bigr) \cdots \bfr\bigl( X_{m_tm_t},
[I\mdd I]
\bigr) = q^{|I\capt M|},
\]  
in view of Lemma \ref{lem2.1}.
The formula for $\bfr\bigl( [I\mdd I], [M\mdd N] \bigr)$ follows via
Lemma \ref{lem1.7}.
\end{proof}


\sectionnew{Initial commutation relations} \label{sec3}

 We now use the
computations of $\bfr(-,-)$ obtained so far to derive some
commutation relations, both to illustrate the method and to doublecheck
the results against known relations in the literature. As in the
previous section, let $i$ and $j$ denote indices in
$\{1,\dots,n\}$, and let $I$, $J$, $M$, $N$ denote index sets
contained in $\{1,\dots,n\}$, with $|I|=|J|$ and $|M|=|N|$.

\begin{subsec} \label{sub3.1}  \textbf{Direct application of
\eqref{eq1.7i}.}  
If we set
$a= X_{ij}$ and $b= [I\mdd J]$ in \eqref{eq1.7i}, we obtain 
\begin{equation} \label{eq3.1}
\sum_{l,L} \bfr\bigl( X_{il}, [I\mdd L] \bigr) X_{lj} [L\mdd J]=
\sum_{l,L} \bfr\bigl( X_{lj}, [L\mdd J] \bigr) [I\mdd L] X_{il}. 
\end{equation} 
We claim that \eqref{eq3.1} reduces to
\begin{equation} \label{eq3.2}
\begin{aligned}
q^{\delta(i,I)}X_{ij} &[I\mdd J] + \bigl( 1-\delta(i,I)
\bigr) \qhat \sum_{\substack{l\in I\\ l<i}} (-q)^{-|(l,i)\capt I|} X_{lj}
[I\setmin l\sqp i\mdd J] =  \\
 &q^{\delta(j,J)}[I\mdd J]X_{ij} + \bigl( 1-\delta(j,J)
\bigr) \qhat \sum_{\substack{l\in J\\ l>j}} (-q)^{-|(j,l)\capt J|} [I\mdd
J\setmin l\sqp j] X_{il}.  
\end{aligned} 
\end{equation}

According to Lemma \ref{lem2.2} and Corollary \ref{cor2.3}, $\bfr\bigl(
X_{il}, [I\mdd L]
\bigr) =0$ unless $i\ge l$ and $I\le L$. By Lemma \ref{lem2.1}, $\bfr\bigl(
X_{ii}, [I\mdd L] \bigr) =0$ unless
$L=I$, and $\bfr\bigl( X_{ii}, [I\mdd I] \bigr)= q^{\delta(i,I)}$. When
$i>l$, Corollary \ref{cor2.5} shows that $\bfr\bigl( X_{il}, [I\mdd L]
\bigr)$ is nonzero only when
$i\in L$, $l\in I$, and $I\setmin l= L\setmin i$. In such cases,
$i\notin I$ and $L= I\setmin l\sqp i$, and the exponent of $-q$ that
appears in \eqref{eq2.7ii} is $-|(l,i)\capt I\capt L|= -|(l,i)\capt I|$.
Thus, the left hand sides of \eqref{eq3.1} and \eqref{eq3.2} agree. 

Similarly, $\bfr\bigl( X_{lj}, [L\mdd J] \bigr) =0$ unless $l\ge j$ and
$L\le J$, while $\bfr\bigl( X_{jj}, [L\mdd J] \bigr) =0$ unless $L=J$, and
$\bfr\bigl( X_{jj}, [J\mdd J] \bigr)= q^{\delta(j,J)}$. When 
$l>j$, Corollary \ref{cor2.5} shows that $\bfr\bigl( X_{lj}, [L\mdd J]
\bigr)$ is nonzero only when $l\in J$, $j\in L\setmin J$, and $L= J\setmin
l\sqp j$. In such cases, the exponent of $-q$ that
appears in \eqref{eq2.7ii} is $-|(j,l)\capt L\capt J|= -|(j,l)\capt J|$.
Therefore, the right hand sides of \eqref{eq3.1} and \eqref{eq3.2} agree.
This establishes \eqref{eq3.2}.
\end{subsec}

\begin{subsec} \label{sub3.2} \textbf{Application of the transpose
automorphism.}  
There are
several ways to obtain a second commutation relation of a similar kind
to \eqref{eq3.2}. First, we could set $a= [I\mdd J]$ and $b= X_{ij}$ in
\eqref{eq1.7i} and proceed as above. Alternatively, we could apply the
automorphism
$\tau$, the anti-automorphism
$\tau_2$, or the isomorphism $\beta$ of \S\ref{sub1.6} to \eqref{eq3.2}
itself. As we shall see in
\S\ref{sub3.4} below, the first three ways are equivalent, up to some
relabelling. The use of $\beta$ is discussed in \S\ref{sub3.5}. 

Among the first three alternatives above, the most
convenient choice is to apply the transpose automorphism
$\tau$ to \eqref{eq3.2}. If we do this, and then relabel the terms by
interchanging $i
\leftrightarrow j$ and $I\leftrightarrow J$, we obtain 
\begin{equation} \label{eq3.3}
\begin{aligned}
q^{\delta(j,J)}X_{ij} &[I\mdd J] + \bigl(
1-\delta(j,J) \bigr) \qhat \sum_{\substack{l\in J\\ l<j}}
(-q)^{-|(l,j)\capt J|} X_{il} [I\mdd J\setmin l\sqp j] =  \\
 &q^{\delta(i,I)}[I\mdd J]X_{ij} + \bigl( 1-\delta(i,I)
\bigr) \qhat \sum_{\substack{l\in I\\ l>i}} (-q)^{-|(i,l)\capt I|}
[I\setmin l\sqp i\mdd J] X_{lj}.  
\end{aligned}   
\end{equation}
\end{subsec}

\begin{subsec} \label{sub3.3} \textbf{Some known cases.}
We now compare some cases of \eqref{eq3.2} and \eqref{eq3.3} with the
literature.

When $i\in I$ and $j\in J$, \eqref{eq3.2} and \eqref{eq3.3} both yield
$qX_{ij} [I\mdd J]= q[I\mdd J]X_{ij}$, the well known fact that $X_{ij}$ and
$[I\mdd J]$ commute in that case. (This is just the centrality of the
quantum determinant in the subalgebra $k\langle X_{st}\mid s\in I,\,
t\in J\rangle$.) If $i\in I$ and $j\notin J$, then \eqref{eq3.2} yields
\begin{equation} \label{eq3.4}
q X_{ij}[I\mdd J] = [I\mdd J]X_{ij} + \qhat \sum_{\substack{l\in J\\ l>j}}
(-q)^{-|(j,l)\capt J|} [I\mdd J\setmin l\sqp j] X_{il}.  
\end{equation}
Multiply \eqref{eq3.4} by $q^{-1}$, and note that $q^{-1}
(-q)^{-|(j,l)\capt J|}= - (-q)^{-|[j,l]\capt J|}$. With this modification,
\eqref{eq3.4} recovers
\cite[Lemma A.1(b)]{GLduke} (which is the second equation of
\cite[Lemma 4.5.1(2)]{PaWa}, rewritten in present notation).
Similarly, consider the case that $i\notin I$ and $j\in J$. Then
\eqref{eq3.3} yields
\begin{equation} \label{eq3.5}
q X_{ij} [I\mdd J] = 
[I\mdd J]X_{ij} + \qhat \sum_{\substack{l\in I\\ l>i}}
(-q)^{-|(i,l)\capt I|} [I\setmin l\sqp i\mdd J] X_{lj}.  
\end{equation} 
We again multiply by $q^{-1}$, and note that $q^{-1} (-q)^{-|(i,l)\capt I|}
= - (-q)^{-|[i,l]\capt I|}$. Thus, \eqref{eq3.5} recovers \cite[Lemma
A.2(c), Equation (A.3)]{GLijm} (which is the second equation of
\cite[Lemma 4.5.1(4)]{PaWa} in present notation). 

Finally, let us consider the case when $i\notin I$ and $j\notin J$. We
may assume that $I\sqp i= J\sqp j= \{1,\dots,n\}$. If we write $\shat=
\{1,\dots,n\}\setmin \{s\}$ for $s=1,\dots,n$, then \eqref{eq3.2} yields
\begin{equation} \label{eq3.6}
X_{ij}[\ihat\mdd \jhat] + \qhat
\sum_{\substack{l\in I\\ l<i}} (-q)^{l+1-i} X_{lj}
[\lhat\mdd \jhat] = [\ihat\mdd \jhat]X_{ij} + \qhat \sum_{\substack{l\in J\\
l>j}} (-q)^{j+1-l} [\ihat\mdd \lhat] X_{il}.   
\end{equation}
Multiplying \eqref{eq3.6} by $q^{-1}$ and then interchanging
$q\leftrightarrow q^{-1}$ recovers the fourth equation of
\cite[Lemma 5.1.2]{PaWa}.
\end{subsec}

\begin{subsec} \label{sub3.4} \textbf{Remark.} 
As mentioned above, \eqref{eq3.3} could also have
been obtained by setting $a= [I\mdd
J]$ and $b= X_{ij}$ in \eqref{eq1.7i} and proceeding as with \eqref{eq3.2}.
In fact, interchanging any choice of $a$ and $b$ in \eqref{eq1.7i} has
the same effect as applying $\tau$, as follows.

First, apply $\tau$ to \eqref{eq1.7i}, and use \eqref{eq1.12} for both $a$
and $b$. This yields
\begin{equation} \label{eq3.7}
\sum_{(a),(b)} \bfr(a_1,b_1) \tau(a)_1 \tau(b)_1= \sum_{(a),(b)}
\bfr(a_2,b_2) \tau(b)_2 \tau(a)_2.  
\end{equation}
Invoking Lemma \ref{lem1.7}, and setting $a'= \tau(a)$ and $b'= \tau(b)$,
\eqref{eq3.7} becomes
\begin{equation} \label{eq3.8}
\sum_{(a'),(b')} \bfr(b'_2,a'_2) a'_1b'_1 = \sum_{(a'),(b')}
\bfr(b'_1,a'_1) b'_2a'_2.  
\end{equation}
Equation \eqref{eq3.8} is nothing but \eqref{eq1.7i} with $a$ and $b$
replaced by
$b'$ and $a'$, respectively.

Similarly, applying the anti-automorphism $\tau_2$ to \eqref{eq1.7i} and
relabelling again recovers \eqref{eq1.7i} with $a$ and $b$ interchanged.
\end{subsec}

\begin{subsec} \label{sub3.5} \textbf{Two further commutation relations.} 
Each case of
commutation relations for $X_{ij}$ and
$[I\mdd J]$ derived in \cite{PaWa} has four subcases -- two pairs in
which one equation of each pair is obtained from the other by inserting
a $q$-Laplace relation. Two commutation relations from each group of
four correspond to our equations \eqref{eq3.2} and \eqref{eq3.3}. It is more
efficient to derive the remaining two by applying the isomorphism
$\beta$ of \S\ref{sub1.6}, as follows. For that purpose, set $A'=
\O_{q^{-1}}(M_n(k))$, and recall the notation $X'_{ij}$ and $[I\mdd
J]'$ for generators and quantum minors in $A'$.

First, consider the relation \eqref{eq3.2} in $A'$, but replace $i$, $j$,
$I$, $J$ by $\itil$, $\jtil$, $\Itil$, $\Jtil$, respectively. The result
is  
\begin{equation} \label{eq3.9}
\begin{aligned} 
q^{-\delta(\itil,\tilde{I})} &X'_{\itil \jtil} [\Itil\mdd
\Jtil]' + \bigl( 1-\delta(\itil,\Itil) \bigr) (-\qhat) \sum_{\substack{\ltil\in \tilde{I}\\
\ltil<\itil}} (-q)^{|(\ltil,\itil)\capt \tilde{I}|} X'_{\ltil\jtil}
[\Itil\setmin \ltil\sqp \itil\mdd \Jtil]' =  \\
 &q^{-\delta(\jtil,\tilde{J})} [\Itil\mdd \Jtil]' X'_{\itil \jtil} +
\bigl( 1-\delta(\jtil,\Jtil) \bigr) (-\qhat) \sum_{\substack{\ltil\in
\tilde{J}\\
\ltil>\jtil}} (-q)^{|(\jtil,\ltil)\capt \tilde{J}|} [\Itil\mdd
\Jtil\setmin \ltil\sqp \jtil]' X'_{\itil \ltil}. 
\end{aligned}  
\end{equation}
Now set 
\begin{alignat}{3}
\itil &= \omega_0(i)  \qquad\qquad\qquad
&\jtil &= \omega_0(j)  \qquad\qquad\qquad
&\ltil &= \omega_0(l)  \notag \\
\Itil &= \omega_0(I)  
&\Jtil &= \omega_0(J)  \notag
\end{alignat}
and apply $\beta^{-1}$ to \eqref{eq3.9}. This yields
\begin{equation} \label{eq3.10}
\begin{aligned}  
q^{-\delta(i,I)} &X_{ij}[I\mdd J] + \bigl( \delta(i,I)-1
\bigr) \qhat \sum_{\substack{l\in I\\ l>i}} (-q)^{|(i,l)\capt I|}
X_{lj}[I\setmin l\sqp i\mdd J] =  \\
 &q^{-\delta(j,J)} [I\mdd J]X_{ij} + \bigl( \delta(j,J)-1 \bigr)
\qhat  \sum_{\substack{l\in J\\ l<j}}
(-q)^{|(l,j)\capt J|} [I\mdd J\setmin l\sqp j] X_{il}.  
\end{aligned}   
\end{equation}

Similarly, the relation \eqref{eq3.3} in $A'$ can be written
\begin{equation} \label{eq3.11}
\begin{aligned} 
q^{-\delta(\jtil,\tilde{J})} &X'_{\itil \jtil} [\Itil\mdd
\Jtil]' + \bigl( 1-\delta(\jtil,\Jtil) \bigr) (-\qhat) 
\sum_{\substack{\ltil\in
\tilde{J}\\ \ltil<\jtil}} (-q)^{|(\ltil,\jtil)\capt \tilde{J}|}
X'_{\itil
\ltil} [\Itil\mdd \Jtil\setmin \ltil\sqp \jtil]' =  \\
 &q^{-\delta(\itil,\tilde{I})} [\Itil\mdd \Jtil]' X'_{\itil \jtil} +
\bigl( 1-\delta(\itil,\Itil) \bigr) (-\qhat) \sum_{\substack{\ltil\in
\tilde{I}\\
\ltil>\itil}} (-q)^{|(\itil,\ltil)\capt \tilde{I}|} [\Itil\setmin
\ltil\sqp \itil\mdd \Jtil]' X'_{\ltil \jtil}. 
\end{aligned}  
\end{equation}
Applying $\beta^{-1}$ to \eqref{eq3.11} as above, we conclude that
\begin{equation} \label{eq3.12}
\begin{aligned} 
q^{-\delta(j,J)} &X_{ij}[I\mdd J] + \bigl(
\delta(j,J)-1 \bigr) \qhat \sum_{\substack{l\in J\\ l>j}}
(-q)^{|(j,l)\capt J|} X_{il} [I\mdd J\setmin l\sqp j] =  \\
 &q^{-\delta(i,I)} [I\mdd J]X_{ij} + \bigl( \delta(i,I)-1
\bigr) \qhat \sum_{\substack{l\in I\\ l<i}} (-q)^{|(l,i)\capt I|}
[I\setmin l\sqp i\mdd J] X_{lj}.  
\end{aligned}  
\end{equation}
\end{subsec}

\begin{subsec} \label{sub3.6} \textbf{Quasicommutation.} 
Elements $a,b\in A$ are said to
{\it quasicommute\/} or {\it $q$-com\-mute\/} provided they commute up to
a power of $q$, that is, $ab= q^mba$ for some integer $m$. The
relations \eqref{eq1.6} say that two of the standard generators for $A$
which have the same row (or column) indices must quasicommute, and it is
natural to expect other instances of this in $A$. From
the results above, we can recover the quasicommutation relations for
quantum minors given by Krob and Leclerc
\cite{KrLe}. These apply to certain quantum minors whose row (or column)
index sets are disjoint. Cases allowing non-disjoint index sets were
obtained by Leclerc and Zelevinsky by investing quantum Pl\"ucker relations
\cite[Lemmas 2.1--2.3]{LeZe}. Building on the results of \cite{LeZe}, 
Scott determined exactly which pairs of
quantum minors quasicommute, and calculated the corresponding relations
\cite[Theorems 1,2]{Sco}. We recover some other cases of his results in
Corollary \ref{newcor5.5} below.

First, consider $X_{ij}$ and $[M\mdd N]$, with $i\in M$.
If $j<\min(N)$, then either \eqref{eq3.3} or \eqref{eq3.10} implies that
$X_{ij}[M\mdd N]= q[M\mdd N]X_{ij}$, while if $j>\max(N)$, then by
either \eqref{eq3.2} or \eqref{eq3.12}, $X_{ij}[M\mdd N]= q^{-1}[M\mdd
N]X_{ij}$. Of course, if $j\in N$, then $X_{ij}[M\mdd N]= [M\mdd N]X_{ij}$.

Now suppose that $I\subseteq M$ and that $J$ and $N$ are {\it separated\/}
in the following sense: there is a partition $J = J'\sqp J''$ such
that
\[
\max(J')< \min(N)\le \max(N)< \min(J'').
\]
Each of the generators $X_{i_{\sigma(l)},j_l}$ occurring in \eqref{eq1.10}
quasicommutes with $[M\mdd N]$ as in the previous paragraph, whence
\[
X_{i_{\sigma(1)},j_1} X_{i_{\sigma(2)},j_2} \cdots
X_{i_{\sigma(t)},j_t} [M\mdd N]= q^{|J'|-|J''|} [M\mdd N]
X_{i_{\sigma(1)},j_1} X_{i_{\sigma(2)},j_2} \cdots
X_{i_{\sigma(t)},j_t}
\]
for all $\sigma\in S_t$. Consequently,
\begin{equation} \label{eq3.13}
[I\mdd J] [M\mdd N]= q^{|J'|-|J''|} [M\mdd N] [I\mdd J] 
\end{equation}
under the present hypotheses. This recovers \cite[Lemma 3.7]{KrLe} (after
interchanging $q$ and $q^{-1}$). In fact, \eqref{eq3.13} holds when
$I\subseteq M$ and $J$ and $N$ are {\it weakly separated\/} in the sense of
\cite{LeZe}, meaning that there is a partition $J\setmin N=
J'\sqp J''$ such that $\max(J')< \min(N\setmin J)\le \max(N\setmin J)<
\min(J'')$ \cite[Lemma 2.1]{LeZe}.

Applying $\tau$ to \eqref{eq3.13} and relabelling, we find that
\begin{equation} \label{eq3.14}
[I\mdd J] [M\mdd N]= q^{|I'|-|I''|} [M\mdd N] [I\mdd J] 
\end{equation}
when $J\subseteq N$ and $I = I'\sqp I''$ with
$\max(I')< \min(M)\le \max(M)< \min(I'')$.
\end{subsec}


\sectionnew{Computation of $\bfr\bigl( [I\mdd J], [M\mdd N] \bigr)$}
\label{sec4}

Throughout this section, let $I$, $J$, $M$, $N$ denote index sets
contained in the interval $\{1,\dots,n\}$, with $|I|=|J|$ and $|M|=|N|$. Our
goal is to develop a formula for $\bfr\bigl( [I\mdd J], [M\mdd N] \bigr)$.

\begin{procl} \label{lem4.1} {\bf Lemma.} 
If $\bfr\bigl( [I\mdd J], [M\mdd N]
\bigr)\ne 0$, then $I\capt M= J\capt N$ and $I\cupt M= J\cupt N$.
\end{procl}

\begin{proof} We induct on $|I|$, starting with
the case $[I\mdd J]= X_{ij}$. If $i=j$, Lemma \ref{lem2.1} implies that
$M=N$, and the conclusion is clear. If $i\ne j$, then $i>j$ by Lemma
\ref{lem2.2}, whence Corollary \ref{cor2.5} implies that $i\in N$, $j\in M$,
and $M\setmin j= N\setmin i$. Consequently, $I\capt M= J\capt N=
\varnothing$ and $I\cupt M= J\cupt N$.

Now suppose that $|I|\ge 2$. If $I=J$, then Proposition \ref{prop2.6}
implies that $M=N$, and we are done. Hence, we may assume that $I\ne J$.
Since
$|I|=|J|$, there must exist an element $j\in J\setmin I$. Set
$J=J_1\sqp J_2$ with
$J_1=\{j\}$ and $J_2=J\setmin j$, and write \eqref{eq1.19} in the
form
\begin{equation} \label{eq4.1}
\pmqdot \bfr\bigl( [I\mdd J], [M\mdd N] \bigr)= \sum_{i\in I} \sum_L
\pmqdot \bfr\bigl( X_{ij}, [M\mdd L] \bigr) \bfr\bigl( [I\setmin i\mdd
J\setmin j], [L\mdd N] \bigr). 
\end{equation}
Since $\bfr\bigl( [I\mdd J], [M\mdd N] \bigr) \ne 0$, \eqref{eq4.1} implies
that
\begin{equation} \label{eq4.2}
\bfr\bigl( X_{ij}, [M\mdd L] \bigr) \bfr\bigl( [I\setmin i\mdd
J\setmin j], [L\mdd N] \bigr)\ne 0  
\end{equation}
for some $i\in I$ and some $L$.

Note that $i\ne j$, because $j\notin I$. Equation \eqref{eq4.2} and Lemma
\ref{lem2.2} now show that
$i>j$, and then Corollary \ref{cor2.5} implies that $i\in L$, $j\in M$, and
$L\setmin i= M\setmin j$. Consequently, $i\notin M$ and $j\notin L$,
while $L= (L\capt M)\sqp i$ and  $M= (L\capt M)\sqp j$. Since the second
factor of \eqref{eq4.2} is nonzero, our induction implies that $(I\setmin
i)\capt L= (J\setmin j)\capt N$ and $(I\setmin i)\cupt L= (J\setmin j)\cupt
N$. Now
\[
I\cupt (L\capt M)= (I\setmin i)\cupt i\cupt (L\capt M)= (I\setmin i)\cupt
L= (J\setmin j)\cupt N,
\]
and so $I\cupt M= I\cupt (L\capt M)\cupt j= J\cupt N$. Since $j\notin
I\cupt L$, we see from the equation $(I\setmin i)\cupt
L= (J\setmin j)\cupt N$ that $j\notin N$. Consequently,
\[
I\capt M= I\capt (M\setmin j)= I\capt (L\setmin i)= (I\setmin i)\capt L=
(J\setmin j)\capt N= J\capt N.
\]
This establishes the induction step.
\end{proof}

\begin{procl} \label{lem4.2} {\bf Lemma.} 
Assume that $I\capt M= J\capt N$ and
$I\cupt M= J\cupt N$.

\textup{(a)} $I\setmin J= N\setmin M$ and $J\setmin
I= M\setmin N$.

\textup{(b)} $\bfr\bigl( [I\mdd J], [M\mdd N] \bigr)= q^{|I\capt M|}
(-q)^{\ell(I;J\capt N)- \ell(J;I\capt M)} \bfr\bigl( [I\setmin
M\mdd J\setmin N], [M\mdd N] \bigr).$
\end{procl} 

\begin{proof} (a) This follows easily from the hypotheses.

(b) Write $J= J_1\sqp J_2$ with $J_1=
J\setmin N$ and  $J_2= J\capt N= I\capt M$, and recall equation
\eqref{eq1.19}. We focus first on the term on the right hand side of
\eqref{eq1.19} with $I_2=J_2$ and $L= N$, in which case
$I_1= I\setmin M$. For this term, we have
\begin{equation} \label{eq4.3}
\begin{aligned} (-q)^{\ell(I_1;I_2)} \bfr\bigl(
[I_1\mdd J_1], [M\mdd L] \bigr)  &\bfr\bigl( [I_2\mdd J_2], [L\mdd N]
\bigr) = \\
 &(-q)^{\ell(I\setmin M; J\capt N)} q^{|I\capt M|} \bfr\bigl( [I\setmin
M\mdd J\setmin N], [M\mdd N] \bigr), 
\end{aligned}  
\end{equation}
in view of Proposition \ref{prop2.6}. We claim that all other terms on the
right hand side of \eqref{eq1.19} vanish.

Suppose that $\bfr\bigl( [I_1\mdd J_1], [M\mdd L] \bigr) 
\bfr\bigl( [I_2\mdd J_2], [L\mdd N] \bigr) \ne 0$ for some $I_1$, $I_2$,
$L$. Lemma \ref{lem4.1} implies that $I_2\capt L= J_2\capt N= J_2$, and then
because
$|I_2|= |J_2|$, we must have $I_2=J_2$. Consequently, Proposition
\ref{prop2.6} implies that $L=N$, verifying the claim. Equations
\eqref{eq1.19} and \eqref{eq4.3} thus yield
\begin{equation} \label{eq4.4}
\begin{aligned}
(-q)^{\ell(J\setmin N;I\capt M)} &\bfr\bigl( [I\mdd J], [M\mdd N] \bigr) =
 \\
 &(-q)^{\ell(I\setmin M; J\capt N)} q^{|I\capt M|} \bfr\bigl( [I\setmin
M\mdd J\setmin N], [M\mdd N] \bigr).  
\end{aligned}   
\end{equation}

Finally, we have
\begin{equation}
\begin{aligned} 
\ell(I;J\capt N) &= \ell(I\setmin M;J\capt N)+ \ell(I\capt M;J\capt
N) \\
\ell(J;I\capt M) &= \ell(J\setmin N;I\capt M)+ \ell(J\capt N;I\capt M),
\end{aligned}  \notag
\end{equation}
and since $I\capt M= J\capt N$, we obtain
\begin{equation} \label{eq4.5}
\ell(I\setmin M;J\capt N)- \ell(J\setmin N;I\capt M)= \ell(I;J\capt N)-
\ell(J;I\capt M).  
\end{equation}
Part (b) follows from \eqref{eq4.4} and \eqref{eq4.5}.
\end{proof}

\begin{procl} \label{lem4.3} {\bf Lemma.} 
Assume that $I\capt M= J\capt N=
\varnothing$ and
$I\cupt M= J\cupt N$. Then
\[
\bfr\bigl( [I\mdd J],[M\mdd N] \bigr)= (-q)^{\ell(I\cupt N; I\setmin J)-
\ell(J\cupt M; J\setmin I)} \bfr\bigl( [I\setmin J\mdd
J\setmin I], [M\setmin N\mdd N\setmin M] \bigr).
\]
\end{procl}

\begin{proof} Write $J= J_1\sqp J_2$ with $J_1=
I\capt J$ and  $J_2= J\setmin I$, and recall \eqref{eq1.19}. Consider the
term with $I_1=J_1$ and $L=M$, in which case $I_2=I\setmin J$. Since
$I_1\capt M=\varnothing$, Proposition \ref{prop2.6} implies that
$\bfr\bigl( [I_1\mdd J_1], [M\mdd L] \bigr) =1$. Thus, for this term of
\eqref{eq1.19}, we have
\begin{equation} \label{eq4.6}
\begin{aligned}
(-q)^{\ell(I_1;I_2)} \bfr\bigl(
[I_1\mdd J_1], [M\mdd L] \bigr)  &\bfr\bigl( [I_2\mdd J_2], [L\mdd N]
\bigr) =  \\
 &(-q)^{\ell(I\capt J; I\setmin J)} \bfr\bigl( [I\setmin J\mdd J\setmin I],
[M\mdd N] \bigr). 
\end{aligned}  
\end{equation}

We next claim that all other terms on the right
hand side of \eqref{eq1.19} vanish. Hence, suppose that $\bfr\bigl(
[I_1\mdd J_1], [M\mdd L] \bigr)  \bfr\bigl( [I_2\mdd J_2], [L\mdd N] \bigr)
\ne 0$ for some $I_1$, $I_2$, $L$. Lemma \ref{lem4.1} implies that $I_2\capt
L= J_2\capt N=
\varnothing$ and $I_2\cupt L= J_2\cupt N= (J\setmin I)\cupt N$, from which
it follows that $I_2= N\setmin L$. Now $I_2\capt J\subseteq N\capt J=
\varnothing$, and so $I_2\subseteq I\setmin J$. Since also
\[
|I_2|= |J_2|= |J\setmin I|= |I\setmin J|,
\]
we must have $I_2= I\setmin J$. Consequently, $I_1=J_1$, and then
Proposition \ref{prop2.6} implies that $L=M$. This verifies the claim. As a
result, \eqref{eq1.19} and \eqref{eq4.6} combine to yield
\begin{equation} \label{eq4.7}
\bfr\bigl( [I\mdd J],[M\mdd N] \bigr)= (-q)^{\ell(I\capt J; I\setmin J)-
\ell(I\capt J; J\setmin I)} \bfr\bigl( [I\setmin J\mdd J\setmin I],
[M\mdd N] \bigr).  
\end{equation}

Note that $(I\setmin J)\capt M= (J\setmin I)\capt N= \varnothing$ and
$(I\setmin J)\cupt M= M\cupt N= (J\setmin I)\cupt N$. Hence, \eqref{eq4.7}
also holds with $I$, $J$, $M$, $N$ replaced by $N$, $M$, $J\setmin I$,
$I\setmin J$, respectively. That is, 
\begin{equation} \label{eq4.8}
\begin{aligned}
\bfr\bigl( [N\mdd M], &[J\setmin I\mdd I\setmin J] \bigr)=  \\
 &(-q)^{\ell(N\capt M; N\setmin M)- \ell(N\capt M; M\setmin N)} \bfr\bigl(
[N\setmin M\mdd M\setmin N], [J\setmin I\mdd I\setmin J] \bigr).  
\end{aligned}  
\end{equation}
In view of Lemma \ref{lem1.7}, \eqref{eq4.8} can be rewritten as
\begin{equation} \label{eq4.9}
\begin{aligned}
\bfr\bigl( [I\setmin J\mdd J\setmin I], &[M\mdd N] \bigr)=  \\
 &(-q)^{\ell(N\capt M; N\setmin M)- \ell(N\capt M; M\setmin N)} \bfr\bigl(
[I\setmin J\mdd J\setmin I], [M\setmin N\mdd N\setmin M] \bigr).  
\end{aligned}  
\end{equation}

Combining \eqref{eq4.7} and \eqref{eq4.9}, we obtain
\begin{equation} \label{eq4.10}
\bfr\bigl( [I\mdd J],[M\mdd N] \bigr)= (-q)^\lambda \bfr\bigl(
[I\setmin J\mdd J\setmin I], [M\setmin N\mdd N\setmin M] \bigr),  
\end{equation}
where (recalling Lemma \ref{lem4.2}(a))
\begin{equation} \label{eq4.11}
\begin{aligned} 
\lambda &= \ell(I\capt J; I\setmin J)-
\ell(I\capt J; J\setmin I)+ \ell(N\capt M; N\setmin M)- \ell(N\capt M;
M\setmin N) \\
 &= \ell((I\capt J)\sqcup (M\capt N); I\setmin J)- \ell((I\capt J)\sqcup
(M\capt N); J\setmin I).
\end{aligned} 
\end{equation}

Next, observe that
\[ 
I\cupt N = (I\setmin J)\sqp (I\capt J)\sqp (M\capt N)  \qquad\qquad
J\cupt M = (J\setmin I)\sqp (I\capt J)\sqp (M\capt N).  
\]
Because $|I\setmin J|= |J\setmin I|$, we have
$\ell(I\setmin J;I\setmin J)= \ell(J\setmin I;J\setmin I)$, and
therefore
\begin{equation} \label{eq4.12}
\lambda= \ell(I\cupt N; I\setmin J)-
\ell(J\cupt M; J\setmin I). 
\end{equation}
Equations \eqref{eq4.10} and \eqref{eq4.12} establish the lemma.
\end{proof}

In view of Lemmas \ref{lem4.1}--\ref{lem4.3}, it only remains to calculate
$\bfr\bigl( [I\mdd J],[M\mdd N] \bigr)$ in case 
\[
 (I\cupt N)\capt (J\cupt M) = \varnothing  \qquad\qquad\qquad\qquad
I\cupt M = J\cupt N, 
\]
 whence $I=N$ and $J=M$.
Further, because of Corollary \ref{cor2.3}, we may assume
that $I> J$. In these cases, certain sums of powers of $-q$ appear in
$\bfr\bigl( [I\mdd J],[M\mdd N] \bigr)$, and we introduce the following
notation to deal with them.

\begin{subsec} \label{sub4.4} \textbf{Definition of $\xi_q(I;J)$.} 
Recall that for
$d\in\NN$, the {\it
$(-q)$-integer\/}
$[d]_{-q}$ is given by
\begin{equation}
\begin{aligned} 
{[d]}_{-q} &= \frac{(-q)^d-(-q)^{-d}}{(-q)-(-q)^{-1}}=
(-q)^{d-1}+(-q)^{d-3} +\cdots+ (-q)^{-(d-1)} \\
 &= (-q)^{1-d}(1+q^2+q^4 +\cdots+ q^{2d-2}). 
\end{aligned} \notag
\end{equation} 
Hence, $1+q^2+q^4 +\cdots+ q^{2d-2}=
(-q)^{d-1}[d]_{-q}$.

We next define a scalar $\xi_q(I;J)$, for index sets $I\ge J$, as
follows. First set $m=|I|$ and write $I= \{r_1< \cdots< r_m\}$. Then
set $d_l= |[1,r_l]\capt J|-l+1$ for $l=1,\dots,m$, noting that $d_l\ge
1$ because $J\le I$. Finally, define
\[
\xi_q(I;J)= [d_1]_{-q} [d_2]_{-q} \cdots [d_m]_{-q},
\]
with the convention that $\xi_q(\varnothing;\varnothing)= 1$.
When $I\capt J=\varnothing$, as in the next lemma, each $d_l=
\ell(r_l;J)-l+1$. Note that $[d]_{-q^{-1}}= [d]_{-q}$ for all
$d\in\NN$, whence $\xi_{q^{-1}}(I;J)= \xi_q(I;J)$.
\end{subsec}

\begin{procl} \label{lem4.5} {\bf Lemma.} 
If $I>J$ and $I\capt J= \varnothing$,
then
\begin{equation} \label{eq4.13}
\bfr\bigl( [I\mdd J],[J\mdd I] \bigr)= \qhat^{|I|} (-q)^{\ell(J;I)-
\ell(I;I)} \xi_q(I;J). 
\end{equation}
\end{procl}

\begin{proof} Set $m=|I|=|J|$, write $I= \{r_1< \cdots< r_m\}$, and set
$d_l= \ell(r_l;J)-l+1$ for $l=1,\dots,m$ as
in \S\ref{sub4.4}.

 We proceed by induction on $m$. If $m=1$, then $J=\{j\}$ for some
$j<r_1$, whence $\ell(J;I)= \ell(I;I)= 0$. Moreover,
$d_1=1$ and so $\xi_q(I;J)=1$. By \eqref{eq1.9}, $\bfr\bigl( [I\mdd
J],[J\mdd I] \bigr)=
\bfr(X_{r_1j}, X_{jr_1})= \qhat$, which verifies \eqref{eq4.13}
in this case.

Now suppose that $m>1$. Write $I= I_1\sqp I_2$ with $I_1=\{r_1\}$ and
$I_2= \{r_2,\dots,r_m\}$. Since $\ell(I_1;I_2)=0$, equation \eqref{eq1.17}
implies that
\[
[I\mdd J]= \sum_{j\in J} (-q)^{\ell(j;J\setmin j)} X_{r_1j} [I_2\mdd
J\setmin j].
\]
Applying \eqref{eq1.7ii}, we obtain
\begin{equation} \label{eq4.14}
\bfr\bigl( [I\mdd J],[J\mdd I] \bigr)= \sum_{j\in J} \sum_L
(-q)^{|[1,j)\capt J|}
\bfr\bigl( X_{r_1j},[J\mdd L] \bigr) \bfr\bigl( [I_2\mdd J\setmin
j],[L\mdd I] \bigr). 
\end{equation} 
According to Lemma \ref{lem2.2} and Corollary
\ref{cor2.5}, a nonzero term can occur on the right hand side of
\eqref{eq4.14} only if
$r_1>j$ and $r_1\in L$, as well as
$J\setmin j= L\setmin r_1$, in which case 
\[
\bfr\bigl( X_{r_1j},[J\mdd L] \bigr)= \qhat (-q)^{|[1,j)\capt
J|- |[1,r_1)\capt L|}.
\]
 Now $|[1,r_1)\capt L|= |[1,r_1)\capt (L\setmin r_1)|=
|[1,r_1)\capt (J\setmin j)|= d_1-1$, 
 and so
\begin{equation} \label{eq4.15}
\bfr\bigl( X_{r_1j},[J\mdd L] \bigr)= \qhat (-q)^{1+ |[1,j)\capt J|-
d_1}. 
\end{equation} 
Next, note that $L= J\setmin j\sqp r_1$, whence $L\capt
I= \{r_1\}$. Consequently, $I_2\capt L=  (J\setmin j)\capt I=
\varnothing$ and $I_2\cupt L = I\cupt L = (J\setmin j)\cupt I$. Lemma
\ref{lem4.3} now implies that 
\begin{equation} \label{eq4.16}
\bfr\bigl( [I_2\mdd J\setmin j],[L\mdd I] \bigr)= (-q)^\lambda
\bfr\bigl( [I_2\mdd J\setmin j], [J\setmin j\mdd I_2] \bigr)
\end{equation}
where
\begin{equation} \label{eq4.17}
\begin{aligned} 
\lambda &= \ell(I;I_2)-\ell(L;J\setmin j) \\
 &= \ell(I_2;I_2)- \ell(J\setmin j;J\setmin j)+ \ell(r_1;I_2)-
\ell(r_1;J\setmin j)
= -d_1+1.  
\end{aligned}  
\end{equation}
Combining equations (4.14)--(4.17), we obtain
\begin{equation} \label{eq4.18}
\bfr\bigl( [I\mdd J],[J\mdd I] \bigr)=  \qhat \sum_{\substack{j\in J\\
j<r_1}} (-q)^{2+2|[1,j)\capt J| -2d_1} \bfr\bigl( [I_2\mdd J\setmin j],
[J\setmin j\mdd I_2] \bigr). 
\end{equation}

It remains to compute $\bfr\bigl( [I_2\mdd J\setmin j], [J\setmin j\mdd
I_2] \bigr)$ for $j\in J$ with $j<r_1$. Observe that
$I_2> J\setmin j$ for any such $j$, so that our induction
hypothesis will apply. Now
\begin{equation}
\begin{aligned} 
\ell(J\setmin j; I_2) &= \ell(J;I_2)= \ell(J;I)-\ell(J;r_1)=
\ell(J;I)-m+d_1 \\
\ell(I_2;I_2) &= \ell(I;I_2)= \ell(I;I) -m+1,  
\end{aligned} \notag
\end{equation}
whence $\ell(J\setmin j;I_2)- \ell(I_2;I_2)=
\ell(J;I)- \ell(I;I) +d_1-1$. For
$l=1,\dots,m-1$, observe that
\[
\ell(r_{l+1};J\setmin j) -l+1= \ell(r_{l+1};J)-l=
d_{l+1},
\]
and consequently $\xi_q(I_2;J\setmin j)= [d_2]_{-q} [d_3]_{-q} \cdots
[d_m]_{-q}$. 
Thus, our induction hypothesis implies that
\begin{equation} \label{eq4.19}
\bfr\bigl( [I_2\mdd J\setmin j], [J\setmin j\mdd
I_2] \bigr)= \qhat^{m-1} (-q)^{\ell(J;I)- \ell(I;I)
+d_1-1} [d_2]_{-q} [d_3]_{-q} \cdots [d_m]_{-q}. 
\end{equation} 
Inserting
\eqref{eq4.19} in \eqref{eq4.18}, we obtain
\begin{equation} \label{eq4.20}
\begin{aligned}
\bfr\bigl( [I\mdd J], &[J\mdd I] \bigr)=  \\
 &\qhat^m (-q)^{\ell(J;I)- \ell(I;I)
+1-d_1} [d_2]_{-q} [d_3]_{-q} \cdots [d_m]_{-q} \sum_{\substack{j\in J\\
j<r_1}} q^{2|[1,j)\capt J|}. 
\end{aligned} 
\end{equation}
The summation appearing in \eqref{eq4.20} is just $\sum_{t=1}^{d_1}
q^{2(t-1)}= (-q)^{d_1-1} [d_1]_{-q}$, whence
\begin{equation} \label{eq4.21}
[d_2]_{-q} [d_3]_{-q} \cdots [d_m]_{-q} \sum_{\substack{j\in J\\ j<r_1}}
q^{2|[1,j)\capt J|}= (-q)^{d_1-1} \xi_q(I;J). 
\end{equation}
Equations \eqref{eq4.20} and \eqref{eq4.21} establish \eqref{eq4.13},
completing the induction step.
\end{proof}

\begin{procl} \label{thm4.6} {\bf Theorem.} 
Let $I,J,M,N\subseteq \{1,\dots,n\}$
with $|I|=|J|$ and $|M|=|N|$.

\textup{(a)} If $\bfr\bigl( [I\mdd J],[M\mdd N] \bigr) \ne 0$, then 
\begin{equation} \label{eq4.22}
I \ge J;  \qquad\qquad
I\capt M = J\capt N;  \qquad\qquad
I\cupt M = J\cupt N. 
\end{equation}

\textup{(b)} If conditions \textup{\eqref{eq4.22}} hold, then
\begin{equation} \label{eq4.23}
\begin{aligned}  
\bfr\bigl( [I\mdd J],[M\mdd N] \bigr) &= q^{|I\capt M|}
\qhat^{|I\setmin J|} (-q)^\lambda \xi_q(I\setmin J;J\setmin I),
\qquad \text{where} \\
\lambda &= \ell\bigl( (J\setmin N)\cupt(M\setmin I);I\setmin J \bigr)-
\ell\bigl( (J\setmin N)\cupt(M\setmin I);J\setmin I \bigr).
\end{aligned}  
\end{equation}
\end{procl}

\begin{proof} (a) Corollary \ref{cor2.3} and Lemma \ref{lem4.1}.

(b) Recall from Lemma \ref{lem4.2} that $I\setmin J= N\setmin M$ and
$J\setmin I= M\setmin N$. If $I=J$, then we must have $M=N$. In
this case,
$\bfr\bigl( [I\mdd J],[M\mdd N] \bigr)= q^{|I\capt M|}$ by Proposition
\ref{prop2.6}, and we are done. Now assume that $I\ne J$, and note that
$I\setmin J> J\setmin I$. We shall need the observations that
\begin{equation}
\begin{aligned}
(I\setmin M)\cupt N &= I\cupt N  \qquad\qquad\qquad
 &(J\setmin N)\cupt M &= J\cupt M \\
(I\setmin M)\setmin(J\setmin N) &= I\setmin J  
 &(J\setmin N)\setmin(I\setmin M) &= J\setmin I.  
\end{aligned} \notag
\end{equation}

 Applying, successively, Lemmas \ref{lem4.2}, \ref{lem4.3}, and
\ref{lem4.5}, we obtain
\begin{equation} \label{eq4.24}
\bfr\bigl( [I\mdd J],[M\mdd N] \bigr)= q^{|I\capt M|}
\qhat^{|I\setmin J|} (-q)^\lambda \xi_q(I\setmin J;J\setmin I), 
\end{equation} 
where
\begin{equation}
\begin{aligned} 
\lambda= \ell(I;J\capt N) &- \ell(J;I\capt M)+ \ell(I\cupt
N;I\setmin J)  \\
 &- \ell( J\cupt M;J\setmin I) + \ell(J\setmin I;I\setmin J)- \ell(I\setmin
J;I\setmin J).
\end{aligned} \notag
\end{equation} 
Observe that $(I\cupt N)\sqp(J\setmin I)= J\cupt N= I\cupt
M= (J\cupt M)\sqp(I\setmin J)$, whence
\begin{equation} \label{eq4.25}
\begin{aligned} 
\ell(I\cupt N;I\setmin J)- \ell( J\cupt M;J\setmin I) +
\ell(J\setmin I; &I\setmin J) - \ell(I\setmin J;I\setmin J)=  \\
 &\ell(J\cupt M;I\setmin J)- \ell(J\cupt M;J\setmin I).
\end{aligned}  
\end{equation}

Next, observe that $I\setmin N= J\setmin M$ and $N\setmin I= M\setmin
J$. Moreover,
\begin{equation}
\begin{aligned}
 I\cupt M &= I\cupt M\cupt N= I\sqp (N\setmin I)\sqp
(M\setmin N) \\
J\cupt N &= J\cupt M\cupt N= J\sqp (M\setmin J)\sqp (N\setmin M), 
\end{aligned} \notag
\end{equation}
and consequently
\begin{equation}
\begin{aligned}
\ell(I;J\capt N) +\ell(N\setmin I;J\capt N) +\ell(M\setmin
N;J\capt N) &= \ell(I\cupt M;J\capt N) \\
\ell(J;I\capt M)+ \ell(M\setmin J;I\capt M) +\ell(N\setmin M;I\capt M) &=
\ell(J\cupt N;I\capt M).
\end{aligned} \notag
\end{equation}
It follows that
\begin{equation} \label{eq4.26}
\begin{aligned}
\ell(I;J\capt N) &- \ell(J;I\capt M)= \ell(N\setmin M;I\capt
M)- \ell(M\setmin N;J\capt N) \\
 &= |N\setmin M|\cdot|I\capt M| - \ell(I\capt M;N\setmin M)  \\
  &\phantom{= |N\setmin M|\cdot|I\capt M|\ }
- |M\setmin N|\cdot|J\capt N|+ \ell(J\capt N;M\setmin N)  \\
 &= \ell(I\capt M;J\setmin I)- \ell(I\capt M;I\setmin J). 
\end{aligned}   
\end{equation} 
Finally, since 
\[
(J\cupt M)\setmin(I\capt M)= \bigl( J\setmin(J\capt N)\bigr) \cupt
\bigl(M\setmin(I\capt M)\bigr)= (J\setmin N)\cupt (M\setmin I),
\]
we conclude from \eqref{eq4.25} and \eqref{eq4.26} that
\begin{equation} \label{eq4.27}
\lambda=  \ell\bigl( (J\setmin N)\cupt(M\setmin I);I\setmin J \bigr)-
\ell\bigl( (J\setmin N)\cupt(M\setmin I);J\setmin I \bigr). 
\end{equation}
In view of \eqref{eq4.24} and \eqref{eq4.27}, the theorem is proved.
\end{proof}

\begin{subsec} \label{newsub4.7} \textbf{Example.}
Let $[I\mdd J]= [45678\mdd 12345]$ and $[M\mdd N]= [123459\mdd 456789]$,
where we have omitted commas between elements of the index sets. It is
clear that $I\ge J$; moreover, $I\capt M= \{4,5\}= J\capt N$ and $I\cupt M=
\{1,\dots,9\}= J\cupt N$. Hence, conditions \eqref{eq4.22} hold. Now
$I\setmin J= \{6,7,8\}$ and $J\setmin I= \{1,2,3\}$, while $(J\setmin
N)\cupt(M\setmin I)= \{1,2,3,9\}$, whence
\[
\ell\bigl( (J\setmin N)\cupt(M\setmin I); I\setmin J \bigr)- \ell\bigl(
(J\setmin N)\cupt(M\setmin I); J\setmin I \bigr)= 3-6= -3.
\]
Since all the elements of $I\setmin J$ are greater than all the elements of
$J\setmin I$, we have
\[
\xi_q(I\setmin J; J\setmin I)= [3]_{-q} [2]_{-q} [1]_{-q}= (q^2+1+q^{-2})
(-q-q^{-1}).
\]
Thus, we conclude from \eqref{eq4.23} that
\[
\bfr\bigl( [I\mdd J],[M\mdd N] \bigr) = q^2 \qhat^3 (-q)^{-3} (q^2+1+q^{-2})
(-q-q^{-1}).
\]
\end{subsec}


\sectionnew{General commutation relations} \label{sec5}

Now that we have formulas for the value of the braiding
form $\bfr$ on pairs of quantum minors, commutation relations follow
readily from property \eqref{eq1.7i}. The following notation for certain
index sets and exponents will be helpful in displaying the results. Recall
the quantities $\ell({-};{-})$ and $\xi_q({-};{-})$ from
\S\ref{newsub1.8} and \S\ref{sub4.4}.

\begin{subsec} \label{sub5.1} \textbf{Definitions of index sets
$\lsind{X}{Y}$ and $\grind{X}{Y}$ and numerical quantities $\L(S,X,Y)$ and
$\L\nat(T,X,Y)$.}  For any
subsets $X$ and $Y$ of
$\{1,\dots,n\}$, define
\begin{equation} \label{eq5.1}
\begin{aligned} 
\lsind{X}{Y} &= \{ S\subseteq X\cupt Y\mid X\capt Y\subseteq
S;\ |S|=|X|;\ S<X \}  \\
\grind{X}{Y} &= \{ T\subseteq X\cupt Y\mid X\capt Y\subseteq
T;\ |T|=|X|;\ T>X \}.  
\end{aligned} 
\end{equation}
In Section \ref{sec6}, we shall need index sets $\lsqind{X}{Y}$ and
$\grqind{X}{Y}$, defined in the same manner.
For any set $S\subseteq X\cupt Y$ such that $X\capt Y\subseteq S$, set
\begin{equation} \label{eq5.2}
S\nat= S_{X,Y}\nat= (X\capt Y)\sqp \bigl( (X\cupt Y)\setmin S
\bigr). 
\end{equation}
Note that if $S\in \lsind{X}{Y}$ or $S\in
\grind{X}{Y}$, then $|S\nat|= |Y|$. Finally, for $S\in \lsind{X}{Y}$
and $T\in \grind{X}{Y}$, define
\begin{equation} \label{eq5.3}
\begin{aligned} 
\L(S,X,Y) &= \ell\bigl( (S\setmin
S\nat)\cupt(Y\setmin X); X\setmin S \bigr)- \ell\bigl( (S\setmin
S\nat)\cupt(Y\setmin X); S\setmin X \bigr)  \\
\L\nat(T,X,Y) &= \ell\bigl( (T\nat\setmin T)\cupt(X\setmin
Y); T\setmin X \bigr)- \ell\bigl( (T\nat\setmin T)\cupt(X\setmin
Y); X\setmin T \bigr).  
\end{aligned} 
\end{equation}

For example, suppose that $X= \{2,3,4,6\}$ and $Y= \{1,3,5\}$. Then
$\lsind{X}{Y}$ consists of those $4$-element subsets $S$ of $\{1,\dots,6\}$
such that $3\in S$ and $S<X$. There are six such sets:
\[
\{1,2,3,4\},\quad \{1,2,3,5\},\quad \{1,2,3,6\},\quad \{1,3,4,5\},\quad
\{1,3,4,6\},\quad \{2,3,4,5\}.
\]
Similarly, $\grind{X}{Y}$ consists of those $4$-element subsets $T$ of
$\{1,\dots,6\}$ such that $3\in T$ and $T>X$. There are two:
$\{3,4,5,6\}$ and $\{2,3,5,6\}$.
Finally, consider the set $S= \{1,2,3,4\} \in \lsind{X}{Y}$. Then $S\nat=
\{3,5,6\}$, and so
\[
\L(S,X,Y) = \ell\bigl( \{1,2,4,5\}; \{6\} \bigr)- \ell\bigl( (\{1,2,4,5\};
\{1\} \bigr)  = 0-3.
\]
\end{subsec}

\begin{procl} \label{thm5.2} {\bf Theorem.} 
Let $I,J,M,N \subseteq \{1,\dots,n\}$
with
$|I|=|J|$ and $|M|=|N|$. Then
\begin{equation} \label{eq5.4}
\begin{aligned}  
q^{|I\cap M|} [I\mdd J][M\mdd N] &+ q^{|I\cap M|} \sum_{S\in
\lsind{I}{M}} \lambda_S [S\mdd J][S\nat\mdd N] =  \\
 &\ q^{|J\cap N|}[M\mdd N][I\mdd J]+ q^{|J\cap N|} \sum_{T\in
\grind{J}{N}} \mu_T [M\mdd T\nat][I\mdd T],  
\end{aligned}  
\end{equation}
where
\begin{equation} \label{eq5.5}
\begin{aligned}  
\lambda_S &= \qhat^{|I\setmin S|} (-q)^{ \L(S,I,M) }
\xi_q(I\setmin S; S\setmin I) \\
\mu_T &= \qhat^{|T\setmin J|} (-q)^{ \L\nat(T,J,N) } \xi_q(T\setmin
J;J\setmin T)  
\end{aligned}  
\end{equation}
for $S\in \lsind{I}{M}$ and $T\in
\grind{J}{N}$.  \end{procl}

\begin{proof} Taking $a=[I\mdd J]$ and $b=[M\mdd N]$ in \eqref{eq1.7i}, we
obtain
\begin{equation} \label{eq5.6}
\sum_{\substack{|S|=|I|\\ |S'|=|M|}} \bfr\bigl( [I\mdd S], [M\mdd S']
\bigr) [S\mdd J][S'\mdd N]=  \sum_{\substack{|T|=|J|\\ |T'|=|N|}}
\bfr\bigl(  [T\mdd J], [T'\mdd N] \bigr) [M\mdd T'][I\mdd T]. 
\end{equation}
In view of
Corollary \ref{cor2.3} and Lemma \ref{lem4.1}, the left hand summation in
\eqref{eq5.6} can be restricted to index sets $S$ and $S'$ such that 
\begin{equation} \label{eq5.7}
\begin{aligned} 
|S| &= |I| &I &\ge S \\
I\capt M &= S\capt S' &\qquad\qquad\qquad\qquad I\cupt M &= S\cupt S'. 
\end{aligned} 
\end{equation}
Proposition \ref{prop2.6} shows that the coefficient of the term with $S=I$
and $S'=M$ is $q^{|I\cap M|}$, and that the terms with $S=I$ and $S'\ne M$
vanish.

The index sets $S$ and $S'$ such that $S\ne I$ and \eqref{eq5.7} hold are
precisely those for which $S\in \lsind{I}{M}$ and $S'= S\nat$. For
these index sets, Theorem \ref{thm4.6} shows that
\[
\bfr\bigl( [I\mdd S], [M\mdd S'] \bigr) = q^{|I\capt M|} \lambda_S.
\]
Thus, the left hand side of \eqref{eq5.6} reduces to the left hand side of
\eqref{eq5.4}.

Similarly, the right hand side of \eqref{eq5.6} reduces to the right hand
side of \eqref{eq5.4}, and the theorem is proved.  \end{proof}

\begin{procl} \label{cor5.4} {\bf Corollary.} 
Let $I,J,M,N \subseteq
\{1,\dots,n\}$ with
$|I|=|J|$ and $|M|=|N|$. Then
\begin{equation} \label{eq5.10}
\begin{aligned}  
q^{|J\capt N|} [I\mdd J][M\mdd N] &+ q^{|J\capt N|}
\sum_{S\in\lsind{J}{N}} \lambda_S [I\mdd S][M\mdd S\nat] = \\
 &\ q^{|I\capt M|} [M\mdd N][I\mdd J]+ q^{|I\capt M|}
\sum_{T\in\grind{I}{M}} \mu_T [T\nat\mdd N][T\mdd J],  
\end{aligned}  
\end{equation}
where
\begin{equation} \label{eq5.11}
\begin{aligned}  
\lambda_S &= \qhat^{|J\setmin S|} (-q)^{ \L(S,J,N) }
\xi_q(J\setmin S; S\setmin J) \\
\mu_T &= \qhat^{|T\setmin I|} (-q)^{ \L\nat(T,I,M) } \xi_q(T\setmin
I;I\setmin T)  
\end{aligned}  
\end{equation}
for $S\in \lsind{J}{N}$ and $T\in
\grind{I}{M}$.  \end{procl}

\begin{proof} Interchange the index sets in the statement of Theorem
\ref{thm5.2} as follows: $I\leftrightarrow J$ and $M\leftrightarrow N$. Then
apply the automorphism $\tau$ to the resulting version of \eqref{eq5.4} to
obtain \eqref{eq5.10} (recall
\eqref{eq1.13}). 

This corollary can also be obtained from Theorem \ref{thm5.2} by
interchanging
$I\leftrightarrow M$ and $J\leftrightarrow N$, in which case one should
also interchange $S\leftrightarrow T\nat$ and $T\leftrightarrow S\nat$.
\end{proof}

\begin{subsec} \label{sub5.3} \textbf{Further quasicommutation.} 
In particular, Theorem
\ref{thm5.2} yields quasicommutation relations of the form $q^{|I\cap M|}
[I\mdd J][M\mdd N] = q^{|J\cap N|}[M\mdd N][I\mdd J]$ in cases where
the index sets
$\lsind{I}{M}$ and $\grind{J}{N}$ are empty. This occurs, for instance, if
either $[I\mdd J]= [1,\dots,r\mdd n{+}1{-}r,\dots,n]$ or $[M\mdd N]= 
[n{+}1{-}r,\dots,n\mdd 1,\dots,r]$, recovering the well known fact that the
northeasternmost and southwesternmost quantum minors are normal elements of
$A$. Moreover,
\begin{equation} \label{eq5.9}
[1,\dots,r\mdd J] [M\mdd 1,\dots,s]= q^{|J\capt [1,s]|- |[1,r]\capt
M|} [M\mdd 1,\dots,s] [1,\dots,r\mdd J], 
\end{equation}
which is part of
\cite[Proposition 1.1]{HoLeSLn} (with $q^2$ replaced by $q$). Also,
\eqref{eq5.9} immediately implies the type A case of \cite[Equation
(10.3)]{BeZe}.

We record the general quasicommutation relations of the above type in the
corollary below. Part (a) recovers one case of \cite[Theorem 2]{Sco}. It
does not seem, however, that the relations
\eqref{eq3.13} and
\eqref{eq3.14} follow directly from
equations such as \eqref{eq5.4} or \eqref{eq5.10}.
\end{subsec} 

\begin{procl} \label{newcor5.5} {\bf Corollary.}
Let $I,J,M,N \subseteq
\{1,\dots,n\}$ with
$|I|=|J|$ and $|M|=|N|$. 

\textup{(a)} If $\max(M\setmin I)< \min(I\setmin M)$ and $\max(J\setmin N)<
\min(N\setmin J)$, then
\begin{equation} \label{neweq5a}
[I\mdd J][M\mdd N]= q^{|I\capt M|- |J\capt N|} [M\mdd N][I\mdd J].
\end{equation}

\textup{(b)} If $\max(I\setmin M)< \min(M\setmin I)$ and $\max(N\setmin J)<
\min(J\setmin N)$, then
\begin{equation} \label{neweq5b}
[I\mdd J][M\mdd N]= q^{|J\capt N|- |I\capt M|} [M\mdd N][I\mdd J].
\end{equation}
\end{procl}

\begin{proof} (a) If $S\in \lsind{J}{N}$, then $S\setmin (J\capt N) <
J\setmin N$, whence
\[
\max\bigl( S\setmin (J\capt N) \bigr) \le \max(J\setmin N)< \min(N\setmin
J).
\]
But then $S$ is disjoint from $N\setmin J$. Since $J\capt N\subseteq
S\subseteq J\cupt N$ and $|S|=|J|$, this forces
$S=J$, which is ruled out by the assumption $S<J$. Thus, $\lsind{J}{N}=
\varnothing$. Similarly, $\grind{I}{M}=\varnothing$, and thus
\eqref{neweq5a} follows from \eqref{eq5.10}.

(b) Interchange $I\longleftrightarrow M$ and $J\longleftrightarrow N$, and
apply part (a).
\end{proof}

\begin{subsec} \label{sub5.5} \textbf{Example.} 
$[n=6]$ Let $J=N=\{1,2,3\}$, and take $I=
\{1,4,5\}$ and $M= \{2,3,6\}$. We first apply Theorem \ref{thm5.2}. Note
that
$\grind{J}{N}$ is empty because $J=N$. For $S\in \lsind{I}{M}$, we
make the following calculations, where commas have been
deleted for the sake of abbreviation (for instance, $\{123\}$ stands for
the index set $\{1,2,3\}$).
\[
\begin{array}{lccccc}
S &\{123\} &\{124\} &\{125\} &\{134\} &\{135\} \\
S\nat &\{456\} &\{356\} &\{346\} &\{256\} &\{246\} \\
I\setmin S &\{45\} &\{5\} &\{4\} &\{5\} &\{4\} \\
S\setmin I &\{23\} &\{2\} &\{2\} &\{3\} &\{3\} \\
(S\setmin S\nat)\cupt(M\setmin I) &\{1236\} &\{12346\} &\{12356\}
&\{12346\} &\{12356\}  \\
\ell\bigl( (S\setmin S\nat)\cupt(M\setmin I); I\setmin S \bigr) &2 &1
&2 &1 &2 \\
\ell\bigl( (S\setmin S\nat)\cupt(M\setmin I); S\setmin I \bigr) &3 &3 &3
&2 &2 \\
\L(S,I,M) &-1&-2&-1&-1&0 \\
\xi_q(I\setmin S;S\setmin I) &{-}q{-}q^{-1} &1 &1 &1 &1 
\end{array}
\]
Consequently, Theorem \ref{thm5.2} implies that
\begin{equation} \label{eq5.12}
\begin{aligned}  
q^3 [236\mdd J] [145\mdd J] &= [145\mdd J] [236\mdd J]
+\qhat^2(-q)^{-1} ({-}q{-}q^{-1}) [123\mdd J] [456\mdd J] \\
 &\qquad +\qhat(-q)^{-2}
[124\mdd J] [356\mdd J] +\qhat(-q)^{-1} [125\mdd J] [346\mdd J] \\
 &\qquad +\qhat(-q)^{-1} [134\mdd J] [256\mdd J] +\qhat [135\mdd J]
[246\mdd J]. 
\end{aligned}   
\end{equation}
The relation \eqref{eq5.12} matches the one calculated by
Fioresi in \cite[Example 2.22]{Fi} (cf.~the first display on page 435,
where one must replace $q$ by $q^{-1}$ to account for the difference
between \eqref{eq1.6} and the relations used in \cite{Fi}).

For contrast, we record the relation obtained from Corollary \ref{cor5.4}
for the current choices of $I$, $J$, $M$, $N$:
\begin{align} \label{eq5.13}
q^3 [145\mdd J] &[236\mdd J] =  \notag \\
 &[236\mdd J] [145\mdd J]
+\qhat [235\mdd J] [146\mdd J] +\qhat(-q)^{-1} [234\mdd J] [156\mdd J] 
\notag \\
 &\qquad +\qhat [136\mdd J] [245\mdd J] +\qhat^2[135\mdd J]
[246\mdd J] +\qhat^2(-q)^{-1} [134\mdd J] [256\mdd J] \\
 &\qquad +\qhat(-q)^{-1} [126\mdd J] [345\mdd J]
+\qhat^2(-q)^{-1} [125\mdd J] [346\mdd J]  \notag \\
 &\qquad +\qhat^2(-q)^{-2} [124\mdd
J] [356\mdd J] +\qhat(-q)^{-4} [123\mdd J] [456\mdd J]. \notag
\end{align}
\end{subsec}

We derive two further relations from Theorem \ref{thm5.2} and Corollary
\ref{cor5.4} with the help of the isomorphism $\beta$ of \S\ref{sub1.6}, as
in \S\ref{sub3.5}. For use in the upcoming proof, note that since $\omega_0$
reverses inequalities of integers, it also reverses the ordering on index
sets: if
$U$ and $V$ are subsets of
$\{1,\dots,n\}$ with $|U|=|V|$, then
$U\le V$ if and only if $\omega_0U\ge \omega_0V$.

\begin{procl} \label{thm5.6} {\bf Theorem.} 
Let $I,J,M,N \subseteq \{1,\dots,n\}$
with
$|I|=|J|$ and $|M|=|N|$. Then
\begin{equation} \label{eq5.14}
\begin{aligned}  
q^{|J\cap N|} [I\mdd J][M\mdd N] &+ q^{|J\cap N|} \sum_{S\in
\grind{I}{M}} \mutil_S [S\mdd J][S\nat\mdd N] =  \\
 &\ q^{|I\cap M|}[M\mdd N][I\mdd J]+ q^{|I\cap M|} \sum_{T\in
\lsind{J}{N}} \lambdatil_T [M\mdd T\nat][I\mdd T],  
\end{aligned}   
\end{equation}
where
\begin{equation} \label{eq5.15}
\begin{aligned}  
\mutil_S &= (-\qhat)^{|S\setmin I|} (-q)^{
-\L\nat(S,I,M) } \xi_q(S\setmin I; I\setmin S) \\
\lambdatil_T &= (-\qhat)^{|J\setmin T|}
(-q)^{ -\L(T,J,N) } \xi_q(J\setmin T; T\setmin J) 
\end{aligned}  
\end{equation}
for $S\in \grind{I}{M}$ and $T\in \lsind{J}{N}$. 
\end{procl}

\begin{proof} Just for this proof, write $\Util= \omega_0U$ for index
sets $U$, and observe that
\begin{equation} 
\begin{aligned}  
\omega_0\bigl( \grind{I}{M} \bigr) &=
\lsind{\Itil}{\Mtil}  &\omega_0\bigl( \lsind{J}{N} \bigr) &=
\grind{\Jtil}{\Ntil}.  
\end{aligned} \notag
\end{equation}
Note also that $\Stil\nat= \widetilde{S\nat}$ for $S\in
\grind{I}{M}$, and similarly $\Ttil\nat= \widetilde{T\nat}$ for
$T\in \lsind{J}{N}$.

Set $A'= \O_{q^{-1}}(M_n(k))$, with generators $X'_{ij}$ and braiding
form $\bfr'$, and label quantum minors in $A'$ in the form $[I\mdd J]'$.
Recall the isomorphism $\beta: A\rightarrow A'$ from \S\ref{sub1.6}, and 
equation \eqref{eq1.14}. Note that when specializing general results to
$A'$, the scalars $q$ and $\qhat$ change to $q^{-1}$ and $-\qhat$,
respectively.

Now apply Theorem \ref{thm5.2} to the quantum minors $[\Itil\mdd\Jtil]'$ and
$[\Mtil\mdd\Ntil]'$ in $A'$. We obtain
\begin{equation} \label{eq5.16}
\begin{aligned}  
q^{-|\Itl\capt\Mtl|} [\Itil\mdd\Jtil]' &[\Mtil\mdd\Ntil]' +
q^{-|\Itl\capt\Mtl|} \sum_{S\in \grind{I}{M}} \lambda'_{\Stl}
[\Stil\mdd\Jtil]' [\Stil\nat\mdd\Ntil]' =  \\
 &\ q^{-|\Jtl\capt\Ntl|} [\Mtil\mdd\Ntil]' [\Itil\mdd\Jtil]'
+q^{-|\Jtl\capt\Ntl|} \sum_{T\in
\lsind{J}{N}} \mu'_{\Ttl} [\Mtil\mdd\Ttil\nat]' [\Itil\mdd\Ttil]', 
\end{aligned} 
\end{equation}
where
\begin{equation} 
\begin{aligned}  
\lambda'_{\Stl} &= (-\qhat)^{|I\setmin S|}
(-q)^{ -\L(\Stl,\Itl,\Mtl) } \xi_q(\Itil\setmin\Stil;
\Stil\setmin\Itil) \\
\mu'_{\Ttl} &= (-\qhat)^{|T\setmin J|} (-q)^{ -\L\nat(\Ttl,\Jtl,\Ntl) }
\xi_q(\Ttil\setmin\Jtil; \Jtil\setmin\Ttil)  
\end{aligned} \notag
\end{equation}
for $S\in \grind{I}{M}$ and $T\in \lsind{J}{N}$. (Here we have
simplified the exponents of the $-\qhat$ terms and invested the
observation that $\xi_{q^{-1}}(U;V)= \xi_q(U;V)$ for any $U$, $V$.)
Applying the isomorphism $\beta^{-1}$ to \eqref{eq5.16} yields
\begin{equation} \label{eq5.17}
\begin{aligned}  
q^{-|I\cap M|} [I\mdd J][M\mdd N] &+ q^{-|I\cap M|}
\sum_{S\in \grind{I}{M}} \lambda'_{\Stl} [S\mdd J][S\nat\mdd N] =  \\
 &\ q^{-|J\cap N|}[M\mdd N][I\mdd J]+ q^{-|J\cap N|} \sum_{T\in
\lsind{J}{N}} \mu'_{\Ttl} [M\mdd T\nat][I\mdd T]  
\end{aligned}  
\end{equation} 
in $A$. Equation \eqref{eq5.14} will follow from \eqref{eq5.17} once we see
that $\lambda'_{\Stl}= \mutil_S$ and $\mu'_{\Ttl}= \lambdatil_T$ for all
$S$ and $T$.

Let $S\in \grind{I}{M}$, and observe that
\begin{equation} \label{eq5.18}
\begin{aligned}  
S\capt S\nat &= I\capt M &\qquad\qquad\qquad\qquad S\cupt
S\nat &= I\cupt M \\
S\nat\setmin M &= I\setmin S &M\setmin S\nat &= S\setmin I. 
\end{aligned}  
\end{equation}
It follows from Theorem \ref{thm4.6} and Lemma \ref{lem1.7} that
\[
q^{-|I\capt M|} \lambda'_{\Stl}= q^{-|\Itl\capt\Mtl|} \lambda'_{\Stl}=
\bfr'\bigl( [\Itil\mdd\Stil]', [\Mtil\mdd\Stil\nat]' \bigr)= \bfr'\bigl(
[M\mdd S\nat]', [I\mdd S]' \bigr).
\]
With the help of \eqref{eq5.18}, a second application of Theorem
\ref{thm4.6} shows that
\[
\bfr'\bigl( [M\mdd S\nat]', [I\mdd S]' \bigr)= q^{-|I\capt M|}
\mutil_S,
\]
 and therefore $\lambda'_{\Stl}= \mutil_S$. Similarly,
$\mu'_{\Ttl}= \lambdatil_T$ for all $T\in  \lsind{J}{N}$, and the
theorem is proved. \end{proof}

The following corollary is obtained from Theorem \ref{thm5.6} in the same
way as Corollary \ref{cor5.4} from Theorem \ref{thm5.2}.

\begin{procl} \label{cor5.7} {\bf Corollary.} 
Let $I,J,M,N \subseteq
\{1,\dots,n\}$ with
$|I|=|J|$ and $|M|=|N|$. Then
\begin{equation} \label{eq5.19}
\begin{aligned}  
q^{|I\cap M|} [I\mdd J][M\mdd N] &+ q^{|I\cap M|} \sum_{S\in
\grind{J}{N}} \mutil_S [I\mdd S][M\mdd S\nat] =  \\
 &\ q^{|J\cap N|}[M\mdd N][I\mdd J]+ q^{|J\cap N|} \sum_{T\in
\lsind{I}{M}} \lambdatil_T [T\nat\mdd N][T\mdd J],  
\end{aligned} 
\end{equation}
where
\begin{equation} \label{eq5.20}
\begin{aligned}  
\mutil_S &= (-\qhat)^{|S\setmin J|} (-q)^{
-\L\nat(S,J,N) } \xi_q(S\setmin J; J\setmin S) \\
\lambdatil_T &= (-\qhat)^{|I\setmin T|}
(-q)^{ -\L(T,I,M) } \xi_q(I\setmin T; T\setmin I) 
\end{aligned}   
\end{equation}
for $S\in \grind{J}{N}$ and $T\in \lsind{I}{M}$. 
\qed\end{procl}


\sectionnew{Some variants} \label{sec6}

Consider the general form of a commutation relation for quantum minors
$[I\mdd J]$ and $[M\mdd N]$, namely an equation that allows a
product $[I\mdd J] [M\mdd N]$ to be replaced by a scalar multiple of
the reverse product $[M\mdd N] [I\mdd J]$, at the cost of some
additional terms. In an equation such as \eqref{eq5.4}, the
additional terms are of two types -- scalar multiples of $[S\mdd
J][S\nat\mdd N]$ and of $[M\mdd T\nat][I\mdd T]$. In some applications, one
type may be more useful than the other. For instance, the {\it prefered
bases\/} constructed in \cite{GLduke} consist of certain products of
quantum minors in which quantum minors with larger index sets must occur to
the left of those with smaller index sets. Thus, if $|I|< |M|$, then
$[M\mdd N] [I\mdd J]$ and the terms $[M\mdd T\nat][I\mdd T]$ are in
preferred order, but $[I\mdd J] [M\mdd N]$ and the terms $[S\mdd
J][S\nat\mdd N]$ are not. A commutation relation in which all the extra
terms are in preferred order can be achieved by iteration -- after a first
application of \eqref{eq5.4}, apply \eqref{eq5.4} to any products $[S\mdd
J][S\nat\mdd N]$ which appear, and continue until all terms have the
desired form. This produces a relation in which $q^{|I\capt M|} [I\mdd
J][M\mdd N]$ is expressed as
$q^{|J\capt N|} [M\mdd N][I\mdd J]$ plus a linear combination of products
$[S\nat\mdd T\nat] [S\mdd T]$ where $S\in \lsqind{I}{M}$ and $T\in
\grqind{J}{N}$. We begin by illustrating the iteration process in
Example \ref{sub6.1} below.

The aim of this section is to derive closed formulas (i.e., without
iterations) for commutation relations of the type just discussed.

\begin{subsec} \label{sub6.1} \textbf{Example.} 
$[n=4]$ Consider $[I\mdd J]= [23\mdd12]$ and
$[M\mdd N]= [14\mdd23]$. First, (5.4) leads to the relation
\begin{equation} \label{eq6.1}
\begin{aligned}  
{[23\mdd12]} &[14\mdd23] -q [14\mdd23][23\mdd12] =  \\
 &q\qhat  [14\mdd12][23\mdd23] - \qhat (-q)^{-1} [12\mdd12][34\mdd23]
-\qhat [13\mdd12][24\mdd23]. 
\end{aligned}  
\end{equation}
The last two terms on the right hand side of \eqref{eq6.1} must now be
treated. Applying \eqref{eq5.4} in each case, we obtain
\begin{align}
[12\mdd12][34\mdd23] &= q [34\mdd23][12\mdd12] + q\qhat 
[34\mdd12][12\mdd23] \tag{6.2)(i}\label{eq6.2i}\\
[13\mdd12][24\mdd23] &= q [24\mdd23][13\mdd12] + q\qhat 
[24\mdd12][13\mdd23] -\qhat [12\mdd12][34\mdd23].
\tag{6.2)(ii}\label{eq6.2ii} 
\end{align}
Note that \eqref{eq6.2ii} contains a term
involving
$[12\mdd12][34\mdd23]$. Hence, we first substitute that equation into
\eqref{eq6.1}, and then combine the two $[12\mdd12][34\mdd23]$ terms, before
substituting \eqref{eq6.2i} into the result. The final
relation is as follows:
\setcounter{equation}{2}
\begin{equation} \label{eq6.3}
\begin{aligned}  
{[23\mdd12]} [14\mdd23] -q [14\mdd23][23\mdd12] &= q\qhat 
[14\mdd12][23\mdd23] -\qhat q [24\mdd23][13\mdd12] \\
 &\qquad -\qhat^2q [24\mdd12][13\mdd23] +\qhat q^2 [34\mdd23][12\mdd12]
\\
 &\qquad +\qhat^2q^2 [34\mdd12][12\mdd23].  
\end{aligned}  
\end{equation}
In each of the terms on the right hand side of \eqref{eq6.3}, the second
factor is of the form $[S\mdd T]$ where $S\in \{23,13,12\}=
\lsqind{I}{M}$ and $T\in \{23,12\}= \grqind{J}{N}$.
\end{subsec}

\begin{procl} \label{lem6.2} {\bf Lemma.} 
Let $s\in \{1,\dots,n-1\}$, and let $B$
and $C$ be the following subalgebras of $A= \OqMn$:
\begin{equation} 
\begin{aligned}  
B &= k\langle X_{ij} \mid 1\le i\le n,\, 1\le j\le s\rangle  \\ 
C &= k\langle X_{ij} \mid 1\le i\le n,\, s+1\le j\le n\rangle. 
\end{aligned}  \notag
\end{equation}
Then the multiplication map
$\mu: B\otimes_k C\rightarrow A$ is a vector space isomorphism.
\end{procl}

\begin{proof} Let $X$, $Y$, and $Z$ be the standard PBW bases of the
respective algebras $B$, $C$, and $A$. Thus,
\begin{equation} 
\begin{aligned}  
X &= \{ (X_{11}^{b_{11}} \cdots X_{1s}^{b_{1s}})
(X_{21}^{b_{21}} \cdots X_{2s}^{b_{2s}}) \cdots (X_{n1}^{b_{n1}} \cdots
X_{ns}^{b_{ns}}) \mid b_{ij}\in \ZZ^+ \} \\
Y &= \{ (X_{1,s+1}^{c_{1,s+1}} \cdots X_{1n}^{c_{1n}})
(X_{2,s+1}^{c_{2,s+1}} \cdots X_{2n}^{c_{2n}}) \cdots
(X_{n,s+1}^{c_{n,s+1}} \cdots X_{nn}^{c_{nn}}) \mid c_{ij}\in \ZZ^+ \}
\\
Z &= \{ (X_{11}^{a_{11}} \cdots X_{1n}^{a_{1n}})
(X_{21}^{a_{21}} \cdots X_{2n}^{a_{2n}}) \cdots (X_{n1}^{a_{n1}} \cdots
X_{nn}^{a_{nn}}) \mid a_{ij}\in \ZZ^+ \}, 
\end{aligned}  \notag
\end{equation}
where the variables occur in each monomial in lexicographic order.
Observe that the monomials $X_{i1}^{b_{i1}} \cdots X_{is}^{b_{is}}$
and $X_{l,s+1}^{c_{l,s+1}} \cdots X_{ln}^{c_{ln}}$ commute whenever
$i>l$. Hence, any product of a monomial from $X$ with a monomial from $Y$ can
be rewritten as follows:
\begin{multline} 
\bigl[ (X_{11}^{b_{11}} \cdots X_{1s}^{b_{1s}})
(X_{21}^{b_{21}} \cdots X_{2s}^{b_{2s}}) \cdots (X_{n1}^{b_{n1}} \cdots
X_{ns}^{b_{ns}}) \bigr] \bigl[ (X_{1,s+1}^{c_{1,s+1}} \cdots
X_{1n}^{c_{1n}}) \cdot  \\
 \shoveright{ (X_{2,s+1}^{c_{2,s+1}} \cdots X_{2n}^{c_{2n}})
\cdots (X_{n,s+1}^{c_{n,s+1}} \cdots X_{nn}^{c_{nn}}) \bigr] } \\
\shoveleft{ \qquad = (X_{11}^{b_{11}} \cdots
X_{1s}^{b_{1s}})(X_{1,s+1}^{c_{1,s+1}}
\cdots X_{1n}^{c_{1n}}) (X_{21}^{b_{21}} \cdots
X_{2s}^{b_{2s}})(X_{2,s+1}^{c_{2,s+1}} \cdots X_{2n}^{c_{2n}}) \cdots
}  \\
 (X_{n1}^{b_{n1}} \cdots X_{ns}^{b_{ns}})(X_{n,s+1}^{c_{n,s+1}}
\cdots X_{nn}^{c_{nn}}).  \notag
\end{multline}
Consequently, $\mu$ maps the set $\{x\otimes y\mid x\in X,\, y\in
Y\}$ bijectively onto $Z$, and the lemma follows. \end{proof}

\begin{procl} \label{thm6.3} {\bf Theorem.}  
Let $I,J,M,N \subseteq \{1,\dots,n\}$
with
$|I|=|J|$ and $|M|=|N|$. Then
\begin{equation} \label{eq6.4}
q^{|I\capt M|} [I\mdd J][M\mdd N] = q^{|J\capt N|} [M\mdd N][I\mdd
J] + q^{|J\capt N|} \sum_{\substack{S\in \lsqind{I}{M}\\ T\in \grqind{J}{N}\\
(S,T)\ne (I,J)}} \lambdatil_S \mu_T [S\nat\mdd T\nat] [S\mdd
T],  
\end{equation}
where
\begin{equation} \label{eq6.5}
\begin{aligned} 
\lambdatil_S &= (-\qhat)^{|I\setmin S|}
(-q)^{ -\L(S,I,M) } \xi_q(I\setmin S; S\setmin I) \\
\mu_T &= \qhat^{|T\setmin J|} (-q)^{ \L\nat(T,J,N) } \xi_q(T\setmin
J;J\setmin T)  
\end{aligned}  
\end{equation}
for $S\in \lsqind{I}{M}$ and $T\in \grqind{J}{N}$. \end{procl}

\noindent\textbf{Remark.} We have isolated the term $q^{|J\capt N|}
[M\mdd N][I\mdd J]$ on the right hand side of \eqref{eq6.4} to emphasize
that this equation is a commutation relation. It may, of course, be
incorporated in the given summation as a term where $(S,T)= (I,J)$, since
$\lambdatil_I \mu_J =1$. 

\begin{proof} Note that the coefficients $\lambda_S$ and $\mu_T$ defined
in \eqref{eq5.5} also depend on $I$, $J$, $M$, $N$. For purposes of the
present proof, we record that dependence by writing
\begin{equation} 
\begin{aligned} 
\lambda^{X,Y}_S &= \qhat^{|X\setmin S|} (-q)^{ \L(S,X,Y) }
\xi_q(X\setmin S; S\setmin X) \\
\mu^{J,N}_T &= \qhat^{|T\setmin J|} (-q)^{ \L\nat(T,J,N) }
\xi_q(T\setmin J;J\setmin T)  
\end{aligned}  \notag 
\end{equation}
for $S\in \lsqind{X}{Y}$ and $T\in
\grqind{J}{N}$. Note that $\lambda_X^{X,Y}=1$ and $\mu_J^{J,N}=1 $.
For $S\in \lsind{I}{M}$, set
\[
\alpha_S^{I,M}= \sum_{\substack{S_1\in \lsind{I}{M}\\ S_2\in
\lsind{S_1}{S_1\nat}\\ \dots\\ S\in \lsind{S_{i-1}}{S_{i-1}\nat}}} (-1)^i \lambda_{S_1}^{I,M} \lambda_{S_2}^{S_1,S_1\nat} \cdots 
\lambda_{S}^{S_{i-1},S_{i-1}\nat},
\]
where we interpret $S_0=I$ and $S_0\nat= M$ in terms where $i=1$.
Finally, set $\alpha_I^{I,M}=1$. We claim that
\begin{equation} \label{eq6.6}
q^{|I\capt M|} [I\mdd J][M\mdd N] = q^{|J\capt N|} \sum_{\substack{S\in
\lsqind{I}{M}\\ T\in \grqind{J}{N}}} \alpha_S^{I,M} \mu^{J,N}_T
[S\nat\mdd T\nat] [S\mdd T].  
\end{equation}

Let $t=|I|$, and let $\N_t$ denote the collection of $t$-element
subsets of $\{1,\dots,n\}$, partially ordered as in \S\ref{sub1.9}. In
proving \eqref{eq6.6}, we proceed by induction on $I$ relative to the
ordering in
$\N_t$. To start, suppose that $I$ is minimal in $\N_t$ (that is, $I=
\{1,\dots,t\}$). In this case, $\lsind{I}{M}$ is empty, and so Theorem
\ref{thm5.2} implies that
\[
q^{|I\cap M|} [I\mdd J][M\mdd N] = q^{|J\cap N|}[M\mdd
N][I\mdd J]+ q^{|J\cap N|} \sum_{T\in
\grind{J}{N}} \mu^{J,N}_T [M\mdd T\nat][I\mdd T],
\]
which verifies \eqref{eq6.6}.

Now suppose that $I$ is not minimal in $\N_t$, but that \eqref{eq6.6} holds
whenever $I$ is replaced by an index set $I'<I$. Theorem \ref{thm5.2}
implies that
\begin{equation} \label{eq6.7}
\begin{aligned} 
q^{|I\cap M|} [I\mdd J][M\mdd N] &= q^{|J\cap N|} \sum_{T\in
\grqind{J}{N}} \mu^{J,N}_T [M\mdd T\nat][I\mdd T] \\
 &\qquad - q^{|I\cap M|}
\sum_{S_1\in \lsind{I}{M}} \lambda^{I,M}_{S_1} [S_1\mdd J] [S_1\nat\mdd
N]. 
\end{aligned} 
\end{equation}
Recall that $S_1\capt S_1\nat= I\capt M$ for $S_1\in \lsind{I}{M}$, by
definition of $S_1\nat$. Hence, our induction hypothesis yields
\begin{equation} \label{eq6.8}
q^{|I\capt M|} [S_1\mdd J][S_1\nat\mdd N] = q^{|J\capt N|}
\sum_{\substack{S\in \lsqind{S_1}{S_1\nat}\\ T\in \grqind{J}{N}}}
\alpha_S^{S_1,S_1\nat}
\mu^{J,N}_T [S\nat\mdd T\nat] [S\mdd T]  
\end{equation}
for all $S_1\in \lsind{I}{M}$. Substitute \eqref{eq6.8} in \eqref{eq6.7},
which yields
\begin{equation} \label{eq6.9}
\begin{aligned} 
q^{|I\cap M|} [I\mdd J][M\mdd N] &= q^{|J\cap N|} \sum_{T\in
\grqind{J}{N}} \mu^{J,N}_T [M\mdd T\nat][I\mdd T] \\
 & \qquad - q^{|J\cap N|}
\sum_{\substack{S_1\in \lsind{I}{M}\\ S\in \lsqind{S_1}{S_1\nat}\\ T\in
\grqind{J}{N}}} \lambda^{I,M}_{S_1} \alpha_S^{S_1,S_1\nat}
\mu^{J,N}_T [S\nat\mdd T\nat] [S\mdd T]. 
\end{aligned}  
\end{equation}
Since $\alpha_I^{I,M}=1$, the coefficients in the first summation of
\eqref{eq6.9} match the corresponding coefficients in \eqref{eq6.6}. The
second summation of \eqref{eq6.9} may be rewritten in the form
\[
q^{|J\cap N|} \sum_{\substack{S\in
\lsind{I}{M}\\ T\in \grqind{J}{N}}} \beta_S \mu^{J,N}_T
[S\nat\mdd T\nat] [S\mdd T],
\]
where each
\[
\beta_S= - \sum_{\substack{S_1\in
\lsind{I}{M}\\ S\in \lsqind{S_1}{S_1\nat}}} \lambda^{I,M}_{S_1}
\alpha_S^{S_1,S_1\nat} = \alpha_S^{I,M}.
\]
Consequently, \eqref{eq6.9} yields \eqref{eq6.6}, establishing the
induction step. This proves \eqref{eq6.6}.

It remains to show that $\alpha_S^{I,M}= \lambdatil_S$ for $S\in
\lsqind{I}{M}$. 

Observe that all quantities appearing in \eqref{eq6.6} involve index sets
contained in the union $I\cupt J\cupt M\cupt N$, and so they remain the
same if we work in $\Oq(M_\nu(k))$ for some $\nu>n$. Hence, there is no
loss of generality in assuming that $n\ge |I|+|M|$. Thus, if we set
\[
J^* = \{n-|I|+1,\dots,n\}  \qquad\qquad\qquad\qquad
N^* = \{1,\dots,|M|\},
\]
we have $\max(N^*) < \min(J^*)$. Note also that $J^*$ is maximal among
$|I|$-element subsets of $\{1,\dots,n\}$. The quantum minors $[U\mdd
N^*]$, for $U \subseteq \{1,\dots,n\}$ with $|U|=|M|$, are homogeneous
elements of distinct degrees with respect to the grading on $A$
discussed in \S\ref{sub1.2}. Hence, the $[U\mdd N^*]$ are linearly
independent over $k$. Similarly, the $[V\mdd J^*]$, for $V\subseteq
\{1,\dots,n\}$ with $|V|=|I|$, are linearly independent, and thus it
follows from Lemma \ref{lem6.2} that the products $[U\mdd N^*][V\mdd J^*]$
are linearly independent over $k$.

Now apply \eqref{eq6.6} to the quantum minors $[I\mdd J^*]$ and $[M\mdd
N^*]$. Since $\grind{J^*}{N^*}$ is empty, we obtain
\begin{equation} \label{eq6.10}
q^{I\capt M|} [I\mdd J^*][M\mdd N^*] = \sum_{S\in \lsqind{I}{M}}
\alpha_S^{I,M} [S\nat\mdd N^*][S\mdd J^*].  
\end{equation}
However, we also have a relation of this type from Corollary \ref{cor5.7},
which may be written in the form
\begin{equation} \label{eq6.11}
q^{I\capt M|} [I\mdd J^*][M\mdd N^*] = \sum_{T\in \lsqind{I}{M}}
\lambdatil_T [T\nat\mdd N^*][T\mdd J^*].  
\end{equation}
Since the products $[S\nat\mdd N^*][S\mdd J^*]$ are linearly
independent, it follows from \eqref{eq6.10} and \eqref{eq6.11} that
$\alpha_S^{I,M}=
\lambdatil_S$ for all $S\in \lsqind{I}{M}$. Therefore \eqref{eq6.6} implies
\eqref{eq6.4}, as desired. \end{proof}

As is easily checked, Theorem \ref{thm6.3} directly yields equation
\eqref{eq6.3}.

We next consider the derivation of new relations from Theorem \ref{thm6.3}.
Unlike the situation in Section \ref{sec4}, however, the methods used there
to prove Corollary \ref{cor5.4} and Theorem \ref{thm5.6} yield the same
result when applied to Theorem \ref{thm6.3}. Hence, we use the method of
Corollary \ref{cor5.4}.

\begin{procl} \label{cor6.4} {\bf Corollary.} 
Let $I,J,M,N \subseteq
\{1,\dots,n\}$ with
$|I|=|J|$ and $|M|=|N|$. Then 
\begin{equation} \label{eq6.12}
q^{|J\capt N|} [I\mdd J][M\mdd N] = q^{|I\capt M|} [M\mdd N][I\mdd
J] + q^{|I\capt M|} \sum_{\substack{S\in \grqind{I}{M}\\ T\in \lsqind{J}{N}\\
(S,T)\ne (I,J)}} \mu_S \lambdatil_T [S\nat\mdd T\nat] [S\mdd
T],  
\end{equation}
where
\begin{equation} \label{eq6.13}
\begin{aligned} 
\mu_S &= \qhat^{|S\setmin I|} (-q)^{ \L\nat(S,I,M) }
\xi_q(S\setmin I;I\setmin S) \\
\lambdatil_T &= (-\qhat)^{|J\setmin T|} (-q)^{ -\L(T,J,N) }
\xi_q(J\setmin T; T\setmin J)  
\end{aligned}  
\end{equation}
for $S\in \grqind{I}{M}$ and $T\in \lsqind{J}{N}$. \end{procl}

\begin{proof} Interchange $I\leftrightarrow J$ and $M\leftrightarrow N$
in the statement of Theorem \ref{thm6.3}, and also interchange the roles of
$S$ and $T$ in the summation. This yields
\begin{equation} \label{eq6.14}
\begin{aligned}
q^{|J\capt N|} [J\mdd I] &[N\mdd M] =  \\
 &q^{|I\capt M|} [N\mdd M][J\mdd
I] + q^{|I\capt M|} \sum_{\substack{T\in \lsqind{J}{N}\\ S\in \grqind{I}{M}\\
(T,S)\ne (J,I)}} \lambdatil^{J,N}_T \mu^{I,M}_S [T\nat\mdd S\nat]
[T\mdd S], 
\end{aligned}  
\end{equation}
where we have placed the superscripts on $\lambdatil^{J,N}_T$ and
$\mu^{I,M}_S$ as reminders of the changes required when carrying over
\eqref{eq6.5} to the present situation. Thus, observe that
$\lambdatil^{J,N}_T$ and
$\mu^{I,M}_S$ are equal to the scalars denoted $\lambdatil_T$ and
$\mu_S$ in \eqref{eq6.13}. Consequently, an application of the automorphism
$\tau$ to \eqref{eq6.14} yields \eqref{eq6.12} (recall \eqref{eq1.13}).
\end{proof}

\begin{subsec} \label{sub6.5} \textbf{Remark.} 
In addition to \eqref{eq6.4} and \eqref{eq6.12}, one can
derive two commutation relations for quantum minors $[I\mdd J]$ and
$[M\mdd N]$ in which the additional terms involve products in the same
order as
$[I\mdd J][M\mdd N]$, rather than in reverse order. To obtain such
results, simply interchange the roles of $[I\mdd J]$ and $[M\mdd N]$ in
Theorem \ref{thm6.3} and Corollary \ref{cor6.4}. One may wish to simplify
the coefficients -- for instance, with the help of observations such as
\eqref{eq5.18}, one sees that $\L(S\nat,M,I)= \L\nat(S,I,M)$. We leave this
to the interested reader.  \end{subsec}

\begin{subsec} \label{sub6.6} \textbf{Example.} 
$[n=4]$ We close the section by applying
Corollary \ref{cor6.4} to the quantum minors $[I\mdd J]= [23\mdd13]$ and
$[M\mdd N]= [14\mdd24]$. In this case, equation \eqref{eq6.12} becomes
\begin{equation} \label{eq6.15}
\begin{aligned} 
{[23\mdd13]} [14\mdd24] &= [14\mdd24][23\mdd13] + \qhat
[13\mdd24][24\mdd13] +\qhat(-q)^{-1} [12\mdd24][34\mdd13] \\
 &\qquad +(-\qhat) [14\mdd34][23\mdd12] +\qhat(-\qhat)
[13\mdd34][24\mdd12] \\
 &\qquad + \qhat(-q)^{-1}(-\qhat) [12\mdd34][34\mdd12]. 
\end{aligned}  
\end{equation}
Equation \eqref{eq6.15} matches the relation calculated by Fioresi in
\cite[Example 6.2]{Fitwo} (after replacing $q$ by $q^{-1}$). 
\end{subsec}


\sectionnew{Poisson brackets} \label{sec7}

In this final section, we use the commutation relations for quantum
minors obtained above to derive expressions for the standard Poisson
bracket on pairs of classical minors in $\OMn$. In particular, we
recover, for the case of the standard bracket, a formula calculated by
Kupershmidt in
\cite{Kup}. Although the study of Poisson brackets is often restricted
to characteristic zero, that restriction is not needed for the results
below.

\begin{subsec} \label{sub7.1} \textbf{Standard Poisson bracket on $\OMn$.} 
Recall that a
{\it Poisson bracket\/} on a commutative $k$-algebra $B$ is a
$k$-bilinear map $\{{-},{-}\}: B\times B\rightarrow B$ such that 
\begin{enumerate}
\item[] $B$ is a Lie algebra with respect to $\{{-},{-}\}$, and
\item[] $\{b,{-}\}$ is a derivation for each $b\in B$.
\end{enumerate}
Note that a Poisson bracket is uniquely determined by its values on
pairs of elements from a $k$-algebra generating set for $B$.

Write $\OMn$ as a
commutative polynomial ring over $k$ in indeterminates $x_{ij}$ for
$i,j = 1,\dots,n$. The {\it standard Poisson bracket\/} on this algebra
is  the unique Poisson bracket such that
\begin{equation} \label{eq7.1}
\begin{aligned} 
\{x_{ij},x_{lj}\} &= x_{ij}x_{lj} &&\qquad (i<l) \\
\{x_{ij},x_{im}\} &= x_{ij}x_{im} &&\qquad (j<m) \\
\{x_{ij},x_{lm}\} &= 0 &&\qquad (i<l,\ j>m) \\
\{x_{ij},x_{lm}\} &= 2 x_{im} x_{lj} &&\qquad (i<l,\ j<m).
\end{aligned}   
\end{equation}
\end{subsec}

\begin{subsec} \label{sub7.2} \textbf{$\Oq(M_n)$ as a quantization of
$\O(M_n)$.}  
It is well
known that $\Oq(M_n(K))$ (for a rational function field $K= k(q)$) is a
quantization of the Poisson algebra
$\bigl( \OMn, \{{-},{-}\} \bigr)$ in the sense that the Poisson
bracket on $\OMn$ is the ``semiclassical limit'' (as $q\rightarrow
1$) of the scaled commutator bracket $\frac1{q-1}[{-},{-}]$ on
$\Oq(M_n(K))$; we indicate the details below.

For the
remainder of this section, replace the scalar $q$ by an indeterminate,
and consider the quantum matrix algebra
$\Oq\bigl( M_n(k(q)) \bigr)$ defined over the rational function field
$k(q)$. The
$k[q^{\pm1}]$-subalgebra $A_0$ of $\Oq\bigl( M_n(k(q)) \bigr)$
generated by the
$X_{ij}$ can be presented (as a $k[q^{\pm1}]$-algebra) by the
generators $X_{ij}$ and relations \eqref{eq1.6}, from which it follows that
there is an isomorphism 
\begin{equation} \label{eq7.2}
A_0/(q-1)A_0 \xrightarrow{\,\cong\,} \OMn  
\end{equation}
sending the cosets $X_{ij}+ (q-1)A_0 \mapsto x_{ij}$ for all $i$, $j$.
We identify $A_0/(q-1)A_0$ with $\OMn$ via \eqref{eq7.2}. Since $\OMn$ is
commutative, the additive commutator $[{-},{-}]$ on $A_0$ takes all its
values in $(q-1)A_0$, and so $\frac1{q-1}[{-},{-}]$ is well-defined on
$A_0$. It follows that the latter bracket induces a well-defined Poisson
bracket on $\OMn$, such that
\begin{equation} \label{eq7.3}
\{\abar,\bbar\}= \overline{ (ab-ba)/(q-1)}  
\end{equation}
for $a,b\in A_0$, where overbars denote cosets modulo $(q-1)A_0$. This
induced bracket is nothing but the standard Poisson bracket on $\OMn$,
as one easily sees by computing its values on pairs of generators
$x_{ij}$,
$x_{lm}$.

We shall apply \eqref{eq7.3} when $\abar$ and $\bbar$ are minors. In order
to reserve the notation $[I\mdd J]$ for classical minors, let us denote
quantum minors in $\Oq\bigl( M_n(k(q)) \bigr)$ in the form $[I\mdd
J]_q$. Note that $[I\mdd J]_q$ is an element of $A_0$, and that the
isomorphism  \eqref{eq7.2} maps the coset of $[I\mdd J]_q$ to $[I\mdd J]$.
Hence, for pairs
of minors, \eqref{eq7.3} can be written as
\begin{equation} \label{eq7.4}
\bigl\{ [I\mdd J], [M\mdd N] \bigr\}= \overline{ \bigl(
[I\mdd J]_q [M\mdd N]_q - [M\mdd N]_q [I\mdd J]_q \bigr) /(q-1)}.
\end{equation}
Combining \eqref{eq7.4} with formulas for additive commutators of quantum
minors thus yields formulas for Poisson brackets of classical minors.
For instance, from \eqref{eq5.9} we obtain
\begin{equation} \label{eq7.5}
\begin{aligned}
\bigl\{ [1,\dots,r\mdd J], &[M\mdd 1,\dots,s] \bigr\} =  \\
 &\bigl( |[1,r]\capt J|- |M\capt [1,s]| \bigr) [1,\dots,r\mdd J] [M\mdd
1,\dots,s],  
\end{aligned}
\end{equation} 
which recovers some cases of \cite[Theorem 2.6]{KoZe}.
\end{subsec}

\begin{procl} \label{thm7.3} {\bf Theorem.} 
Let $I,J,M,N \subseteq \{1,\dots,n\}$
with
$|I|=|J|$ and $|M|=|N|$. Then
\begin{equation} \label{eq7.6}
\begin{aligned}
\bigl\{  [I\mdd J], [M\mdd N] \bigr\} &= \bigl( |J\capt N|-
|I\capt M| \bigr) [I\mdd J] [M\mdd N] \\
 &\qquad +2 \sum_{\substack{j\in J\setmin N\\ n\in N\setmin J\\ j<n}} (-1)^{|(J\dlt N)\capt (j,n)|} [I\mdd J\sqp n\setmin j] [M\mdd
N\sqp j\setmin n] \\
 &\qquad -2 \sum_{\substack{i\in I\setmin M\\ m\in M\setmin I\\ i>m}} (-1)^{|(I\dlt M)\capt (m,i)|} [I\sqp m\setmin i\mdd J][M\sqp
i\setmin m\mdd N].  
\end{aligned}  
\end{equation}
\end{procl}

\begin{proof} Write \eqref{eq5.4} in the form
\begin{equation} \label{eq7.7}
\begin{aligned}
{[I\mdd J]}_q [M\mdd N]_q - [M\mdd N]_q [I\mdd J]_q  &= \bigl(
q^{|J\cap N|- |I\cap M|} -1 \bigr) [M\mdd N]_q [I\mdd J]_q  \\
 &\qquad + q^{|J\cap N|- |I\cap M|}\sum_{T\in
\grind{J}{N}} \mu_T [M\mdd T\nat]_q [I\mdd T]_q  \\
 &\qquad - \sum_{S\in \lsind{I}{M}} \lambda_S [S\mdd J]_q [S\nat\mdd
N]_q . 
\end{aligned}  
\end{equation}
Since $\qhat^2/(q-1)$ vanishes modulo $q-1$, we only need to consider
the terms in the sums for $T\in
\grind{J}{N}$ with $|T\setmin J|=1$ and $S\in \lsind{I}{M}$ with
$|I\setmin S|=1$. Any such $T$ has the form $T= J\sqp n\setmin j$ with
$j\in J\setmin N$ and $n\in N\setmin J$ such that $j<n$, whence $T\nat=
N\sqp j\setmin n$ and $(T\nat\setmin T)\cupt (J\setmin N)= (J\dlt
N)\setmin n$ and so
\begin{equation} 
\begin{aligned}
\L\nat(T,J,N) &= \ell\bigl( (J\dlt N)\setmin n; n\bigr) -
\ell\bigl( (J\dlt N)\setmin n; j \bigr) \\
 &= \ell(J\dlt N;n)- \ell(J\dlt
N;j)+1= -|(J\dlt N)\capt (j,n)|. 
\end{aligned}  \notag
\end{equation}
Similarly, the indices $S$ that appear have the form $S= I\sqp
m\setmin i$ with $i\in I\setmin M$ and $m\in M\setmin I$ such that
$i>m$, whence $S\nat= M\sqp i\setmin m$ and $\L(S,I,M)= -|(I\dlt
M)\capt (m,i)|$. Consequently, dividing \eqref{eq7.7} by $q-1$ and then
reducing the resulting equation modulo $q-1$ yields \eqref{eq7.6}.
\end{proof}

Similarly, Corollary \ref{cor5.4} yields the following result.

\begin{procl} \label{thm7.4} {\bf Theorem.} 
Let $I,J,M,N \subseteq \{1,\dots,n\}$
with
$|I|=|J|$ and $|M|=|N|$. Then
\begin{equation} \label{eq7.8}
\begin{aligned}
\bigl\{  [I\mdd J], [M\mdd N] \bigr\} &= \bigl( |I\capt M|-
|J\capt N| \bigr) [I\mdd J] [M\mdd N] \\
 &\qquad +2 \sum_{\substack{i\in I\setmin M\\ m\in M\setmin I\\ i<m}} (-1)^{|(I\dlt M)\capt (i,m)|} [I\sqp m\setmin i\mdd J][M\sqp
i\setmin m\mdd N] \\
 &\qquad -2 \sum_{\substack{j\in J\setmin N\\ n\in N\setmin J\\ j>n}} (-1)^{|(J\dlt N)\capt (n,j)|} [I\mdd J\sqp n\setmin j] [M\mdd
N\sqp j\setmin n]. \qquad\square 
\end{aligned} 
\end{equation}
\end{procl}

Finally, provided $k$ does not have characteristic 2, we can average
equations \eqref{eq7.6} and \eqref{eq7.8} to obtain the equation
below. 

\begin{procl} \label{cor7.5} {\bf Corollary.} 
Let $I,J,M,N \subseteq
\{1,\dots,n\}$ with
$|I|=|J|$ and $|M|=|N|$. If $\operatorname{char}(k)\ne 2$, then
\begin{equation} \label{eq7.9}
\begin{aligned}
\bigl\{  [I\mdd J], [M\mdd N] \bigr\} &= \sum_{\substack{i\in
I\setmin M\\ m\in M\setmin I\\ i<m}} (-1)^{|(I\dlt M)\capt (i,m)|}
[I\sqp m\setmin i\mdd J][M\sqp i\setmin m\mdd N] \\
 &\qquad - \sum_{\substack{i\in I\setmin M\\ m\in M\setmin I\\ i>m}} (-1)^{|(I\dlt M)\capt (m,i)|} [I\sqp m\setmin i\mdd J][M\sqp
i\setmin m\mdd N] \\
 &\qquad + \sum_{\substack{j\in J\setmin N\\ n\in N\setmin J\\ j<n}} (-1)^{|(J\dlt N)\capt (j,n)|} [I\mdd J\sqp n\setmin j] [M\mdd
N\sqp j\setmin n] \\
 &\qquad - \sum_{\substack{j\in J\setmin N\\ n\in N\setmin J\\ j>n}} (-1)^{|(J\dlt N)\capt (n,j)|} [I\mdd J\sqp n\setmin j] [M\mdd
N\sqp j\setmin n]. \qquad\square 
\end{aligned}  
\end{equation}
\end{procl}

Equation \eqref{eq7.9} is the standard case of Kupershmidt's formula
\cite[Equation (9)]{Kup}. To obtain the standard Poisson bracket in his
setting, make the following choices for the structure constants:
\[
r^{ij}_{lm}= \begin{cases}
 \ 1 &\quad (i>j,\, l=j,\, m=i) \\ 
 -1 &\quad (i<j,\, l=j,\, m=i) \\ 
 \ 0 &\quad \text{(otherwise)}. \end{cases}
\]


\section*{Acknowledgement}

We thank T. H. Lenagan, L. Rigal, A. Zelevinsky, and the referee for their
comments and suggestions concerning this project.


\end{document}